\documentclass[sn-mathphys,Numbered]{sn-jnl}

\usepackage{graphicx}
\usepackage{multirow}
\usepackage{amsmath,amssymb,amsfonts}
\usepackage{amsthm}
\usepackage{mathrsfs}
\usepackage[title]{appendix}
\usepackage{xcolor}
\usepackage{textcomp}
\usepackage{manyfoot}
\usepackage{booktabs}
\usepackage{algorithm}
\usepackage{algorithmicx}
\usepackage{algpseudocode}
\usepackage{listings}
\usepackage{graphicx}
\usepackage{dsfont}
\usepackage{tikz}
\usepackage{mathtools}

\newtheorem{theorem}{Theorem}
\newtheorem{prop}[theorem]{Proposition}

\newtheorem{rem}{Remark}

\newtheorem{defi}{Definition}

\newtheorem{coro}{Corollary}[section]
\newtheorem{lem}{Lemma}[section]

\newcommand{\BE}{\begin{equation}}
\newcommand{\EE}{\end{equation}}

\newcommand{\BEN}{\begin{eqnarray}}
\newcommand{\EEN}{\end{eqnarray}}
\newcommand{\BENE}{\begin{eqnarray*}}
\newcommand{\EENE}{\end{eqnarray*}}

\newcommand{\TPA}{T_{\alpha}^p (x_0)}
\newcommand{\HFP}{h_{f}^p (x_0)}
\newcommand{\EFP}{E_f^{(p)}}

\newcommand{\RR}{{\mathbb R}}

\newcommand{\ZZ}{{\mathbb Z}}
\newcommand{\NN}{{\mathbb N}}
\newcommand{\PP}{{\mathbb P}}
\newcommand{\CCC}{{\mathcal C}}
\newcommand{\ep}{\varepsilon}
\newcommand{\Card}{\mbox{Card}}

\newcommand{\BD}{\medskip\begin{defi}}
\newcommand{\ED}{\end{defi}\medskip}
\newcommand{\BR}{\medskip\begin{rem}}
\newcommand{\ER}{\end{rem}\medskip}
\newcommand{\BC}{\medskip\begin{coro}}
\newcommand{\EC}{\end{coro}\medskip}
\newcommand{\BL}{\medskip\begin{lem}}
\newcommand{\EL}{\end{lem}\medskip}
\newcommand{\BP}{\medskip\begin{prop}}
\newcommand{\EP}{\end{prop}\medskip}
\newcommand{\BT}{\medskip\begin{theorem}}
\newcommand{\ET}{\end{theorem}\medskip}

\begin{document}

\title[Characterization of $p$-exponents by continuous wavelet transforms, applications to the multifractal analysis of sum of random pulses]{Characterization of $p$-exponents by continuous wavelet transforms, applications to the multifractal analysis of sum of random pulses}

%%=============================================================%%
%% Prefix	-> \pfx{Dr}
%% GivenName	-> \fnm{Joergen W.}
%% Particle	-> \spfx{van der} -> surname prefix
%% FamilyName	-> \sur{Ploeg}
%% Suffix	-> \sfx{IV}
%% NatureName	-> \tanm{Poet Laureate} -> Title after name
%% Degrees	-> \dgr{MSc, PhD}
%% \author*[1,2]{\pfx{Dr} \fnm{Joergen W.} \spfx{van der} \sur{Ploeg} \sfx{IV} \tanm{Poet Laureate} 
%%                 \dgr{MSc, PhD}}\email{iauthor@gmail.com}
%%=============================================================%%

\author*[1]{\fnm{Guillaume} \sur{Saës}}\email{guillaume.saes@u-pec.fr}

\affil*[1]{\orgdiv{Laboratoire}, \orgname{d'Analyse et de Mathématiques Appliquées CNRS, UMR 8050, UPEC}, \orgaddress{\street{61 Av. du Général de Gaulle}, \city{Créteil}, \postcode{94000}, \country{France}}}

%\affil[2]{\orgdiv{Department}, \orgname{Organization}, \orgaddress{\street{Street}, \city{City}, \postcode{10587}, \state{State}, \country{Country}}}

\abstract{The theory of orthonormal wavelet bases is a useful tool in multifractal analysis, as it provides a characterization of the different exponents of pointwise regularities (Hölder, $p$-exponent, lacunarity, oscillation, etc.). However, for some homogeneous self-similar processes, such as sums of random pulses (sums of regular, well-localized functions whose expansions and translations are random), it is easier to estimate the spectrum using continuous wavelet transforms. In this article, we present a new characterization of $p$-exponents by continuous wavelet transforms and we provide an application to the regularity analysis of sums of random pulses.}


\maketitle

\section{Introduction}\label{sec1}

In the mid-1980s, measurements of the speed of a turbulent flow \cite{Frisch85, Gagne87} highlighted  signals that seemed  irregular everywhere and whose irregularities seemed to vary from one point to another. The notion of pointwise regularity allows to characterize  this type of behavior. For a function $f\in L_{\text{loc}}^{\infty} (\RR)$ the pointwise regularity at a point $x_0$ is most often measured via the Hölder pointwise exponent defined as follows.

\BD
Let $f\in L_{\text{loc}}^{\infty} (\RR)$. Let $x_0\in\RR$ and $\alpha\geq 0$. A function $f$ belongs to $C^\alpha (x_0)$ when there exists a polynomial $P_{f,x_0}$ of degree less than $\alpha$ and $C,r>0$ such that
$$
\forall x \in (x_0-r,x_0+r), \quad |f(x)-P_{f,x_0}(x-x_0)| \leq C |x-x_0|^{\alpha}.
$$
The Hölder exponent of $f$ at $x_0$ is $h_f (x_0)=\sup\{\alpha \ : \ f\in C^{\alpha} (x_0)\}$.
\ED

The Hölder exponent is well defined only when $f$ is locally bounded. Calderón and Zygmund introduced a notion of pointwise regularity for functions that belong to $L_{\text{loc}}^{p}(\RR)$, $p\geq 1$, to study the regularity of solutions of certain partial differential equations \cite{Calderon61}.

\BD\label{defp}
Let $f\in L_{\text{loc}}^{p} (\RR)$ with $p\geq 1$. Let $x_0\in\RR$. A function $f$ belongs to $\TPA$ when there exist a polynomial $P_{f,x_0}$ of degree less than $\alpha$ and $C,R>0$ such that
\BE\label{defTp}
\forall r\in (0,R), \quad \left( \frac{1}{r} \int_{x_0-r}^{x_0+r} |f(x)-P_{f,x_0}(x-x_0)|^p dx \right)^{\frac{1}{p}} \leq C r^{\alpha}. \nonumber
\EE
The $p$-exponent of $f$ at $x_0$ is $\HFP =\sup\{\alpha : f\in \TPA\}$.
\ED

\BR If $f\in L^p(\RR)$ then for all $x_0\in\RR$, $f\in T_{-\frac{1}{p}}^p (x_0)$. So for all functions $f\in L^p (\RR)$, we will be interested only in knowing whether $f\in\TPA$ for $\alpha\geq -1/p$.
\ER

The polynomial $P_{f,x_0}$ is unique and note also that, if $\alpha \leq 0$, the polynomial is null. The condition for a function to belong to $\TPA$ extends to values $p \in (0,1)$, in this case the spaces $L^p$ are replaced by Hardy spaces $H^p$ \cite{Jaffard06a}. Note that when $p=+\infty$, then $h^p_f (x_0)=h_f(x_0)$.

The wavelet coefficients are known to be convenient tools to compute Hölder exponents and $p$-exponents \cite{Jaffard04b}. However, for some homogeneous self-similar processes, estimating  some continuous wavelet transforms (with respect to its wavelet coefficients)  turns out to be easier. 

In this paper, our goal is to determine a characterization of  Hölder pointwise exponents and  $p$-exponents by the continuous wavelet transforms (CWT) and to apply it to determine the pointwise regularity of a class of stochastic processes. For $C^\alpha (x_0)$ spaces, a characterizations of the Hölder exponent  was proposed  in \cite{Jaffard06a, Holschneider91} while for $\TPA$ spaces and $p$-exponents comparable results remain to be established. Such characterizations provide  new tools to compute the multifractal $p$-spectrum of certain functions, processes or signals.

\BD
Let $f\in L_{\text{loc}}^{p} (\RR)$ with $p\geq 1$. The multifractal $p$-spectrum $D^{(p)}_f : [-1/p,+\infty] \rightarrow [0,1]\cup \{-\infty\}$ of $f$ is the mapping
\BE
D_f^{(p)} (h) = \dim_H \EFP (h), \nonumber
\EE
where the  set $\EFP (h)$ is defined as
\BE
\EFP (h) = \{x\in\RR : \quad h_f^{(p)} (x)=h\}, \nonumber
\EE
and $\dim_H$ is the Hausdorff dimension with the convention $\dim_H (\emptyset) = -\infty$.
\ED
In the case $p=+\infty$, since $h^{+\infty}_f (x_0)=h_f(x_0)$ we recover the classical definition of the multifractal spectrum.

\BD
Let $f\in L_{\text{loc}}^{\infty} (\RR)$. The multifractal spectrum $D_f : [0,+\infty] \rightarrow \RR_+\cup \{-\infty\}$ of $f$ is the mapping defined for every $h\geq 0$ by
$$
D_f (h) = \dim_H E_f (h),
$$
where the set $E_f (h)$ is defined as $E_f (h) = \{x\in\RR : \quad h_f (x)=h\}.$
\ED

A classic example of a multifractal function is the Riemann series defined by 
$$
R(x) = \sum _{n=1}^{+\infty} \frac{\sin (\pi n^2 x)}{n^2}, \quad x\in\RR.
$$
The multifractal spectrum associated with the Hölder exponents was studied by S. Jaffard \cite{Jaffard96a} and the author already made good use of the continuous wavelet transforms to do so. The spectrum he obtained is the following:
$$
D_{R} (h) = \left\{ 
\begin{tabular}{l l}
     $4h-2$ & if $h\in [\frac{1}{2}, \frac{3}{4}]$ \\
     \\
     $0$ & if $h=\frac{3}{2}$\\
     \\
     $-\infty$ & else. 
\end{tabular}
\right.
$$

On the other hand, the only known result related to the $p$-spectrum of the Riemann function concerns its extension
$$
F_s (x) = \sum_{n=1}^{+\infty} \frac{e^{2 i \pi n^2 x}}{n^s}, \quad x\in\RR.
$$
In \cite{Seuret17}, S. Seuret and A. Ubis  managed to estimate a part of the multifractal $2$-spectrum for $s\in (1/2,1]$. More precisely, they proved that
$$
D_{F_s}^{(2)} (h) = 4h+2-2s, \qquad h\in \left[0, \frac{s}{2}-\frac{1}{4}\right],$$
which leaves open the conjecture that the $p$-spectrum of these functions is independent of $p$ and shifts by $s$ when we operate a fractional derivation of order 1 opened.

Another example of not locally bounded multifractal function is the Brjuno function. Its complex version  was introduced in 1971 by A. Brjuno to give a condition for the function to be holomorphic in $0$. The Brjuno function is defined for every $x\in\mathbb{R}\setminus \mathbb{Q}$ by
$$
B(x) = \sum_{n=0}^{+\infty} xA(x)\ldots A^{n-1} (x) \log \left( \frac{1}{A^n (x)}\right),
$$
where $A$ is the Gaussian map which to an irrational $x$ of $(0,1]$ associates $A(x)=\{1/x\}$ the fractional part of $1/x$. Moreover, S. Marmi, P. Moussa and J-C. Yoccoz \text{ showed in } \cite{Marmi97a}  that, $B$ is not locally bounded and  $B\in L^p (\RR)$ for all $p<+\infty$ (and is even BMO). The multifractal analysis of this function (using $p$-exponents) was achieved by S. Jaffard and B. Martin in \cite{Jaffard18}, and they obtained
$$
D_{B}^{(1)} (h) = \left\{ 
\begin{tabular}{l l}
     $2h$ & if $h\in [0,\frac{1}{2}]$ \\
     $-\infty$ & else. 
\end{tabular}
\right.
$$
As in the case of the Riemann function, the $p$-spectrum of the Bjruno function turns out to be independent of $p$. Note that a continuous wavelet transform technique was already used in the mentioned paper to obtain a lower-bound for the $p$-exponents of $B$.

\medskip

In the two previous cases, the proofs rely on continuous wavelet transform estimates since orthonormal wavelet bases, based on a dyadic grid, are not suitable for the problem. 

\medskip

Continuous wavelet transforms are used to estimate the pointwise regularity of some processes \cite{Saes20}, functions \cite{Jaffard96a, Jaffard18} or signals \cite{Arneodo08a}. In the case of sums of random pulses introduced and studied in many papers (see \cite{Cioczek95a, Cioczek95b, Cioczek96a, Lovejoy85a, Mandelbrot95a}) in order to model rain fields and the volume of water over time \cite{Lovejoy85a}. In these different cases, the properties of almost sure convergence, continuity, stationarity and self-affinity have been studied. During his PhD., Y. Demichel \cite{Demichel06} computed the uniform Hölder exponent of a large family of sums of pulses and was interested in obtaining some information about the structure of the graphs of such functions and in particular, in estimating their Hausdorff dimension. For the study of the pointwise regularity of such processes, a first multifractal analysis has been performed in \cite{Saes20, Saes21}, which will  be the application case we will consider. This last model is defined as a sum of random dilatations and translations of a function arbitrarily chosen at start.

\BD
Let $(\Omega,\mathcal{F},\PP)$ a probability space. Let $(C_n)_{n\in\NN^*}$ be a real Poisson point process whose intensity is the Lebesgue measure on $\RR_+$. Let $S$ be an independent point process of $(C_n)_{n\in\NN^*}$ whose intensity is the Lebesgue measure on $\RR_+^*\times [0,1]$. We write $S=(B_n,X_n)_{n\in\NN^*}$ where $(B_n)_{n\in\NN^*}$ is an increasing sequence.
\ED
By construction, the three sequences of random variables $(C_n)_{n\in\NN^*}$, $(B_n)_{n\in\NN^*}$ and $(X_n)_{n\in\NN^*}$ are independent. Let us now recall the definition the sums of random pulses \cite{Saes20} that we will later study.

\BD\label{defpulse}
Let $\psi : \RR \rightarrow \RR$ a non-zero Lipschitz function with support equal to $[-1,1]$. The sum of random pulses  $F_{\alpha,\eta} : \RR \rightarrow\RR$ is the stochastic process defined by
\BE
\label{defP}
F_{\alpha,\eta}(x) = \sum_{n=1}^{+\infty} C_n^{-\alpha} \psi\left(B_n^{\frac{1}{\eta}}(x-X_n)\right), \ \ \ x\in\RR
\EE
\ED

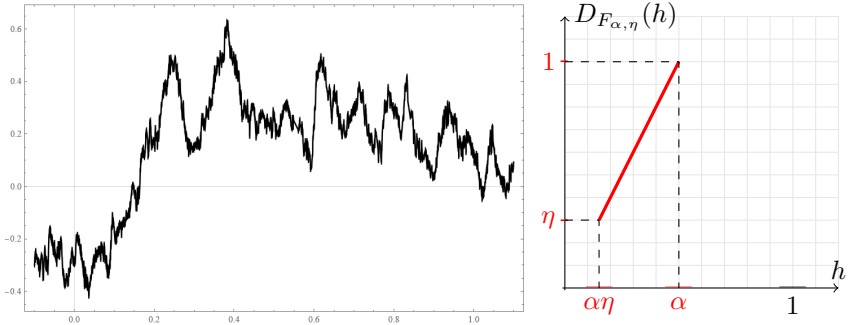
\begin{figure}
    \centering
    \includegraphics[scale=0.28]{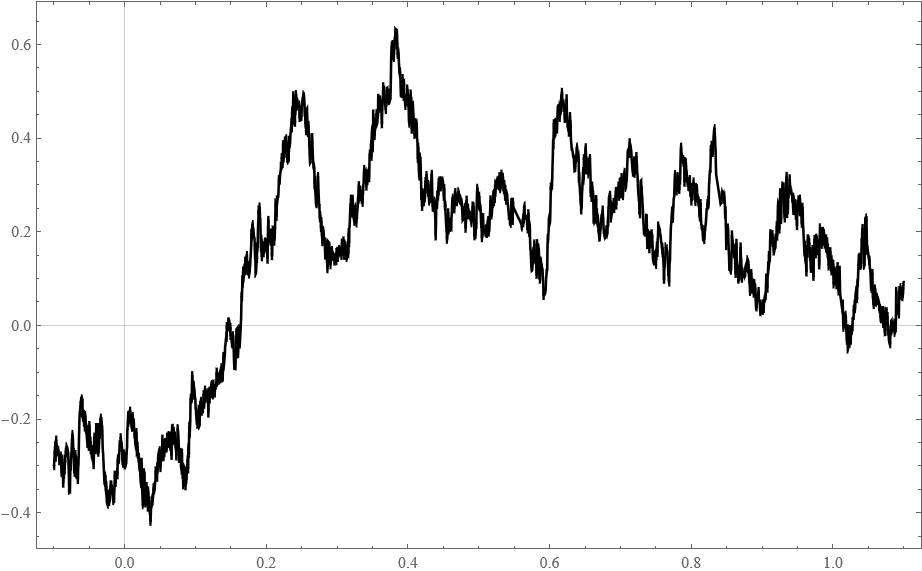}
    \begin{tikzpicture}[scale=0.3]

    \draw[-,very thin,gray!20] (1,0) -- (1,12);
    \draw[-,very thin,gray!20] (2,0) -- (2,12);
    \draw[-,very thin,gray!20] (3,0) -- (3,12);
    \draw[-,very thin,gray!20] (4,0) -- (4,12);
    \draw[-,very thin,gray!20] (5,0) -- (5,12);
    \draw[-,very thin,gray!20] (6,0) -- (6,12);
    \draw[-,very thin,gray!20] (7,0) -- (7,12);
    \draw[-,very thin,gray!20] (8,0) -- (8,12);
    \draw[-,very thin,gray!20] (9,0) -- (9,12);
    \draw[-,very thin,gray!20] (10,0) -- (10,12);
    \draw[-,very thin,gray!20] (11,0) -- (11,12);
    \draw[-,very thin,gray!20] (12,0) -- (12,12);
    \draw[-,very thin,gray!20] (0,1) -- (12,1);
    \draw[-,very thin,gray!20] (0,2) -- (12,2);
    \draw[-,very thin,gray!20] (0,3) -- (12,3);
    \draw[-,very thin,gray!20] (0,4) -- (12,4);
    \draw[-,very thin,gray!20] (0,5) -- (12,5);
    \draw[-,very thin,gray!20] (0,6) -- (12,6);
    \draw[-,very thin,gray!20] (0,7) -- (12,7);
    \draw[-,very thin,gray!20] (0,8) -- (12,8);
    \draw[-,very thin,gray!20] (0,9) -- (12,9);
    \draw[-,very thin,gray!20] (0,10) -- (12,10);
    \draw[-,very thin,gray!20] (0,11) -- (12,11);
    \draw[-,very thin,gray!20] (0,12) -- (12,12);
    \draw[->,black] (-0.1,0) -- (12,0)  node[above] {$h$};
    \draw[->,black] (0,-0.1) -- (0,12)  node[right] {$D_{F_{\alpha,\eta}} (h)$};
    
    \node[red] at (0,3) {-};
    \node[red] at (0,10) {-};
    \draw[red] (0,3) node[left] {$\eta$};
    \draw[red] (0,10) node[left] {$1$};
    \node[red] at (1.5,0) {|};
    \node[red] at (5,0) {|};
    \draw[red] (1.5,0) node[below] {$\alpha\eta$};
    \draw[red] (5,0) node[below] {$\alpha$};

    \draw[] (10,0) node[below] {$1$};
    \node[] at (10,0) {|};
    
    \draw[-,red, very thick] (1.5,3) -- (5,10);

    \draw[dashed] (0,3) -- (1.5,3) -- (1.5,0);
    \draw[dashed] (0,10) -- (5,10) -- (5,0);
    
    \end{tikzpicture}
    \caption{Sample path of $F_{\alpha,\eta}$ with $\alpha=0.5$, $\eta=0.9$ and $\psi:t\mapsto t(1-t^2)^2$ if $t\in [-1,1]$, $0$ else (left) and multifractal spectrum of $F_{\alpha,\eta}$ with $\alpha=0.5$ and $\eta=0.9$ (right)}
    \label{fig:spect1}
\end{figure}

If $\alpha>0$ and $\eta\in (0,1)$, the sample paths of $F_{\alpha,\eta}$ are locally bounded and we represent on Figure \ref{fig:spect1} an example of such a trajectory.

Unlike in the case of random series or lacunar wavelets series \cite{Aubry02,Jaffard00}, the dyadic network is not privileged in the case of  sums of random pulses and it is then not surprising that continuous wavelet transforms are more suited to estimate the multifractal properties in such cases. We emphasize  that the proof of the following result obtained in \cite{Saes20} already partially relied on estimates of some continuous wavelet transforms.

\BT 
Let $\psi:\RR\rightarrow\RR$ be a non-zero Lipschitz function  supported on $[-1,1]$, $\alpha,\eta\in (0,1)$ and let $F_{\alpha,\eta}$ be the random series defined by \eqref{defP}. With probability $1$,
$$
D_{F_{\alpha,\eta}} (h) = \left\{ 
\begin{tabular}{l l}
     $\frac{h}{\alpha}$ & if $h\in [\alpha\eta, \alpha]$ \\
     $-\infty$ & else. 
\end{tabular}
\right.
$$
\ET
The multifractal spectrum is plotted in Figure \ref{fig:spect1}. 

The proof relies on a characterization of the Hölder exponent by continuous wavelet transforms. As it allows to choose wavelets that are ``positioned'' at any point, this tool is more flexible and convenient to use in the case of sums of random pulses  as we can choose a wavelet suited to the location of our pulses. The main difficulty is to obtain a characterization of the $p$-exponents by  continuous wavelet transforms. Section \ref{sec3}  presents a consistent definition of continuous $p$-leaders in terms of wavelet coefficients and  Theorem \ref{propcar}, which together with Theorem \ref{pthm1} is the main result of this paper, yields  a characterization of the  $\TPA$ spaces based on $p$-leaders. Section \ref{sec4} proposes a proof of this theorem.

\medskip

The only known examples of random processes for which one gets almost surely a  $p$-spectrum are lacunary wavelet series \cite{Jaffard15a}.  We will study the $p$-spectrum of $F_{\alpha,\eta}$ processes for $\alpha <0$. Those are examples of processes whose sample paths are almost surely non-locally bounded. As such they are an interesting model since many signals cannot be modeled by locally bounded functions \cite{Jaffard16, Leonarduzzi16} and such situations are in general not well understood yet. The section \ref{sec5} of the article proposes a proof of the following theorem 

\BT\label{pthm1} Let $\eta\in (0,1)$, $\alpha <0$ such that $\eta -1 < \alpha \eta$ and $p\in (1,-1/(\alpha\eta) + 1/\alpha)$. 

Let $\psi:\RR\rightarrow\RR$ a continuous lipschitzian function with support on $[-1,1]$. Let $F_{\alpha,\eta}:\RR\rightarrow\RR$ the process defined by \eqref{defP}. Then almost certainly, we have
$$D_{F_{\alpha,\eta}}^{(p)} (H) = \dfrac{H\eta p+\eta}{\alpha\eta p+1}\ \  \mbox{ if } H\in \left[\alpha\eta, \alpha + \frac{1-\eta}{\eta p} \right].$$
\ET

Finally, we emphasize also  that the $p$-exponents can differ from the Hölder exponents and  thus give additional information on the nature of the irregularities at a point, even in the case of processes with locally bounded sample paths.
It is then natural to study new characterizations of $\TPA$ spaces by continuous wavelet transforms. This was achieved by V. Perrier and C. Basdevant in their paper \cite{Perrier96} in the case of the space $L^p$ (using a result of Stein's book \cite{Stein93}).
We will start from this characterization to obtain the one of $L^p$ spaces as it has been done in the discrete case in the proof of S. Jaffard and C. Melot in \cite{Jaffard05}.

\section{The continuous \texorpdfstring{$p$}{p}-leaders, characterisation of \texorpdfstring{$\TPA$}{TPA}}\label{sec3}

Let $\alpha\in\RR$ be fixed and $r\geq \lfloor\max(1,\alpha) \rfloor$. Let  $\psi\in C^{r}(\RR)$ be a non-zero even function which is called wavelet, with support included in $[-1,1]$, and having $r+1$ vanishing moments, i.e. $\int_\RR \psi (x) dx = \int_\RR x\psi (x) dx =\dots = \int_\RR x^r \psi (x) dx =0$. 

\BD
The continuous wavelet transform of a function $f\in L^2 (\RR)$ is defined in \cite{Murenzi08} (see also \cite{Meyer87, Daubechies92}) by
\BE\label{eqcomp}
W_f (a,b) = \frac{1}{a} \int_{\RR} f(x) \psi\left( \frac{x-b}{a} \right) dx. \nonumber
\EE
\ED

There exists a constant $c_\psi>0$ such that the reconstruction formula is valid for $f\in L^2 (\RR)$ with 
$$
    \ f(x)=\frac{1}{c_\psi} \int_{\RR_+^*}\int_{\RR}\frac{1}{a^{2}} W_f (a,b) \psi\left( \frac{x-b}{a}\right) dbda,\ x\in\RR^d,
$$
if the wavelet $\psi$ verifies the following admissibility condition
\BE\label{AdFour} 
\int_{\RR_{-}} \frac{ | \hat{\psi} (\xi)|^2 }{|\xi|}  d\xi = \int_{\RR_{+}} \frac{ | \hat{\psi} (\xi)|^2 }{|\xi|}  d\xi < + \infty, \nonumber
\EE
which is the case for $\psi$ under our assumption. Additionally it is possible to take another wavelet $\phi$ of class $C^N (\RR)$ where $N>\max(0,\alpha)$ with this time (at least) the first vanishing moment such that the recomposition of $f$ is
\BE\label{eqrerecomp}
f(x) = \int_{\RR_+^*}\int_{\RR}\frac{1}{a^2} W_f(a,b) \phi\left( \frac{x-b}{a}\right)dbda,\ x\in\RR,
\EE
(see for example \cite{Grossmann02, Daubechies92, Jaffard96b} for possible choices of $\phi$). The reconstruction formula \eqref{eqrerecomp} holds almost everywhere if $f \in L^p (\RR)$ with $p >1$. The choice of a reconstruction wavelet that differs from the analyzing wavelet offers an additional flexibility that will prove important in proving Theorem \ref{propcar}.

Before stating the characterization of $\TPA$ spaces by continuous wavelet transforms, we recall the definitions of wavelet coefficients and its associated quantities in an orthonormal wavelet basis.

\BD 
Let $\varphi$ be a non-zero oscillating function, with support in $[-1,1]$, having a number $r_{\varphi}\geq 1$ of vanishing moments and of class $C^{r_{\varphi}-1}(\RR)$. The function $\varphi$ is called a mother wavelet when $\{\varphi_{j,k}(t)= 2^{\frac{j}{2}}\varphi(2^{j}t-k)\}_{(j,k)\in \NN\times\ZZ}$ forms an orthonormal basis of $L^2 (\RR)$. The discrete wavelet coefficients of a function $f\in L^2 (\RR)$ are defined by
$$ c_{j,k}=2^{\frac{j}{2}}\int_{\RR} f(x) \varphi_{j,k}(x)dx=2^{j} \int_{\RR} f(x) \varphi \left (2^j x-k \right)dx, \qquad (j,k)\in\NN\times\ZZ.$$
\ED

For any $(j,k)\in\NN\times\ZZ$, let $\lambda_{j,k}$  be the dyadic interval defined by $\lambda_{j,k}=[k 2^{-j}, (k+1)2^{-j}]$. We will use the notation $3\lambda_{j,k}$ for the union of $\lambda_{j,k}$ and the $2$ intervals adjacent to $\lambda_{j,k}$, $3\lambda_{j,k} =[(k-1)2^{-j}, (k+2)2^{-j}]$.
Finally, for all $j\in\NN$ and $x_0\in\RR$, let $\lambda_j (x_0)$ the unique dyadic interval $\lambda_{j,k}$ such that $x_0\in\lambda_{j,k}$. The wavelet leaders are defined as follows.

\BD
Let $p\in\RR_+^*$ and $f\in L^{p}_{loc} (\RR)$. If $p=+\infty$, then the wavelet leaders of $f$ are
$$\forall (j,k)\in\NN\times\ZZ, \ l_{f}(j,k)=\sup_{\substack{\lambda_{j',k'}\subseteq 3\lambda_{j,k} \\ j' \geq j}} \{|c_{j',k'}|\}.$$
If $p<+\infty$, then the $p$-leaders of $f$ are
$$
\forall (j,k)\in\NN\times\ZZ,\ l^{(p)}_f (j,k) = \left( \sum_{\substack{\lambda_{j',k'}\subseteq 3\lambda_{j,k} \\ j' \geq j}} |c_{j',k'}|^p 2^{-(j'-j)} \right)^{\frac{1}{p}}.   
$$
\ED

In the case of continuous wavelet transforms, we define a continuous version of the $p$-leaders as a local $L^p$-norm of coefficients $W_f (a,b)$. For this, we rely on a characterization of $L^p$-spaces for $p>1$ that V. Perrier and C. Basdevant \cite{Perrier96} (Theorem 3.1) have constructed.

\BT\label{caractth}
Let $f\in L_{\text{loc}}^p (\RR)$ with $p\in (1,+\infty)$. Let $\alpha > -1/p$.

Let $\psi$ be even function with support in $[-1,1]$ having $r\geq \max (\alpha, 1)$ vanishing moments and $\psi\in C^{r} (\RR)$.
There exists $C_1,C_2>0$, depending only on the wavelet $\psi$ such that if 
\[N_f = \left( \int_{-\infty}^{+\infty} \left( \int_0^{+\infty} |W_f (s,t)|^2 \frac{ds}{s} \right)^{\frac{p}{2}} dt \right)^{\frac{1}{p}},\]
then
\[C_2 N_f \leq \lVert f \rVert_{L^p} \leq C_1 N_f.\]
\ET

\BD\label{defpl}
Let $f\in L^{\infty}_{loc} (\RR)$. The continuous leaders of $f$ are
\BE\label{deflead}
  \forall (a,b)\in\RR_+^* \times \RR, \ L_f (a,b) = \sup_{(s,t)\in (0,a) \times B(b,a)} |W_f(s,t)|. \nonumber
\EE
where $B(b,a)=(b-a,b+a)$.

Let $f\in L^p_{loc} (\RR)$. The continuous $p$-leaders of $f$ are
\BE\label{defplead}
  \forall (a,b)\in\RR_+^* \times \RR, \ L_f^{(p)} (a,b) = \left(\frac{1}{a} \int_{B(b,a)} \left(\int_{0}^a |W_f (s,t)|^2 \frac{ds}{s} \right)^{\frac{p}{2}} dt \right)^{\frac{1}{p}}. \nonumber 
\EE
\ED

This indicates that the notion of $\TPA$ regularity can be related to $p$-leaders in the framework of continuous wavelet transforms. However for $p=1$, one cannot obtain such a characterization since, for $p>1$ their characterization is based on the characterization of $L^p$-spaces, and $L^1$ cannot be characterized by wavelets ($L^1$ has no unconditional basis) \cite{Meyer90}. The result we will prove is the following.

\BT\label{propcar}
Let $f\in L_{\text{loc}}^p (\RR)$ with $p\in (1,+\infty)$. Let $\alpha > -1/p$ and $x_0\in\RR$.

Let $\psi$ be even function with support in $[-1,1]$ having $r\geq \max (\alpha, 1)$ vanishing moments and $\psi\in C^{r} (\RR)$.
\begin{enumerate}
    \item For any $f\in\TPA$, there exists $C_f>0$ such that for any $a\in\RR_+^*$ small enough,
\BE\label{hypL}
    |L^{(p)}_f (a,x_0)| \leq C_f  a^\alpha .
\EE
\item Conversely, if $\alpha\notin\NN$ and \eqref{hypL} is verified, then $f$ belongs to $\TPA$.
\end{enumerate}
\ET

Moreover, this theorem allows to determine the $p$-exponent in the following way.

\BC
Let $f\in L_{\text{loc}}^p (\RR)$ with $p\in (1,+\infty)$. If $h^{p}_{x_0} (x_0) \leq 1$ then the $p$-exponent of $f$ in $x_0$ is
\[h_f^{(p)}(x_0) = \liminf_{a\rightarrow 0^+} \frac{\log (|L_f^{(p)} (a,x_0)|)}{\log(a)}.\]
\EC

The section \ref{sec4} present a proof of theorem \ref{propcar}.

\section{Proof of characterization of \texorpdfstring{$\TPA$}{TPA} spaces}\label{sec4}

The Definition \ref{defpl} of $p$-leaders is a continuous version of the discrete leaders defined using wavelet coefficients. We will check on an example the consistency of this definition with classical results on discrete $p$-leaders. 

Consider the cusp function $g_\alpha : x\mapsto |x|^\alpha$ with $\alpha>0$. When $\alpha$ is not an even integer and if we take the wavelet $\varphi$ supported in $[-1,1]$, then the wavelet coefficients in $(j,k)$ are given by a change of variable $u=2^j x-k$ by
\BE 
c_{j,k} = 2^j \int_\RR |x|^\alpha \varphi( 2^j x-k)dx =2^{-\alpha j} \int_{-1}^{1} |u+k|^\alpha \varphi(u) du = 2^{-\alpha j} w_{\varphi,\alpha}(k) \nonumber
\EE
where $w_{\varphi,\alpha} (t) = \int_{-1}^{1} |u+t|^\alpha \varphi(u)du$. We deduce the $p$-leaders of $g_\alpha$ by
\BENE
[l_{g_\alpha}^{(p)} (j,0)]^p &=& 2^j \sum_{\substack{\lambda_{j',k'}\subseteq 3\lambda_{j,0} \\ j' \geq j}} 2^{-(1+\alpha p) j'} w_{\varphi,\alpha} (k')^p  \\
& = & 2^{-\alpha p j} \sum_{\substack{\lambda_{j',k'}\subseteq 3\lambda_{j,0} \\ j' \geq j}} 2^{-(1+\alpha p) (j'-j)} w_{\varphi,\alpha}(k')^p \\
& = & C_\alpha^p 2^{-\alpha j p}
\EENE
where 
\BE 
C_\alpha = \left( \sum_{\substack{\lambda_{j'',k'}\subseteq 3\lambda_{j,0} \\ j'' \geq 0}} 2^{-(1+\alpha p) j''} w_{\varphi,\alpha} (k')^p\right)^{\frac{1}{p}} >0 \nonumber
\EE
which is a well defined constant because $\alpha\geq -1/p$. As a conclusion, it is known \cite{Jaffard89, Jaffard16} that the discrete $p$-leader $l_{g_\alpha}^{(p)} (j,k)$ have a scaling law behavior $2^{-\alpha j}$ in $0$, i.e., there exists $C_\alpha>0$ such that for any $j\in\NN$, $l_{g_\alpha}^{(p)}(j,0) = C_\alpha 2^{-\alpha j}$  and from the discrete characterization of $\TPA$ \cite{Jaffard16}, we have $h_{g_\alpha}^{(p)} (0)=\alpha$. 

In the continuous case of $p$-leaders, it is therefore expected that their behavior follows the same $a^\alpha$ scaling invariance. The continuous wavelet transform of the cusp for $(a,b)\in \RR_+^* \times \RR$ is
$$W_{g_\alpha} (a,b) =\int_{\RR} g_{\alpha}(x) \psi\left( \frac{x-b}{a}\right)\frac{dx}{a} =\int_{\RR}g_\alpha (au+b)\psi(u)du = a^\alpha \int_{-1}^{1} \left|u+\frac{b}{a}\right|^\alpha \psi(u) du.$$
Thus, $W_{g_\alpha}(a,b)= a^\alpha w_{\psi,\alpha} (b/a)$ où $w_{\psi,\alpha} (t)=\int_{-1}^{1} |u+t|^\alpha \psi(u)du$.
The continuous $p$-leaders of the cusp $g_\alpha$ are computed by changing the variable $t=t'/ad$ and $s=s'/a$ in the following manner 
\BENE
L_{g_\alpha}^{(p)} (a,0)^p &=& \frac{1}{a} \int_{-a}^{a} \left(\int_{0}^a | s^\alpha w_{\psi,\alpha} (t/s)|^2 \frac{ds}{s} \right)^{\frac{p}{2}} dt \\
&=& a^{\alpha p} \int_{-1}^{1} \left( \int_0^1 \left| s'^\alpha w_{\psi,\alpha} \left( \frac{t'}{s'} \right)\right|^2 \frac{ds'}{s'}\right)^{\frac{p}{2}}dt'.
\EENE
We obtain precisely the power law of the continuous $p$-leader with $L_{g_\alpha}^{(p)} (a,0) = K a^{\alpha}$ where 
$$K=\left(\int_{-1}^{1} \left( \int_0^1 \left| s'^\alpha w_{\psi,\alpha} \left( \frac{t'}{s'} \right)\right|^2 \frac{ds'}{s'}\right)^{\frac{p}{2}}dt\right)^{\frac{1}{p}}.$$
By integration by part, one easily checks that $w$ is a regular and well-localized function so we deduce the same behavior as the $p$-leaders in the discrete case.

In a first step, we present a proof of part 1. of Theorem \ref{propcar} and in a second step the proof of part 2.

\subsection{Proof of point 1. of Theorem \ref{propcar}}\label{subsec4.1}

We recall the hypothesis of Theorem \ref{propcar}. Let $f\in L_{\text{loc}}^p (\RR)$ with $p\in (1,+\infty)$. Let $\alpha > -1/p$ and $x_0\in\RR$. Let $\psi$ be even function with support in $[-1,1]$ having $r\geq \max (\alpha, 1)$ vanishing moments and $\psi\in C^{r} (\RR)$.

In this section, we prove for any $f\in\TPA$, there exists $C_f>0$ such that for any $a\in\RR_+^*$ small enough, we have \eqref{hypL} holds.

\begin{proof}
Since $f\in \TPA$, there exist two constants $C,R\in\RR_+^*$ and a polynomial $P_{f,x_0}$ of degree less than $\alpha$ such that \eqref{defTp} holds. Recall that if $\alpha \in [-1/p,0]$, then the polynomial $P_{f,x_0}$ vanishes. Let $r\in (0,R)$ and $g$ be the function defined by
\BE\label{eqg}
    g(x)=[f(x)-P_{f,x_0} (x-x_0)]\mathds{1}_{B(x_0,r)}(x), \ x\in\RR. \nonumber
\EE
Thus from \eqref{defTp}, 
\BE\label{eqnorm}
\lVert g \rVert_{L^p} = \left(\int_{x_0-r}^{x_0+r} |f(x)-P_{f,x_0}(x-x_0)|^p dx\right)^{\frac{1}{p}}\leq C r^{\alpha+\frac{1}{p}}.
\EE
According the Theorem \ref{propcar}, there exist $C_1,C_2>0$ depending only on the wavelet $\psi$ such that, 
\BE\label{caracLp}
   C_2 N_g \leq \lVert g \rVert_{L^p} \leq C_1 N_g.
\EE

Let $r>0$. Since the wavelet $\psi$ is supported in $[-1,1]$, $x\mapsto \psi \left( \frac{x-t}{s}\right)$ is supported in $B(t,s)$ and for all $s \in \left[0,\frac{r}{2}\right], t\in B\left(x_0,\frac{r}{2}\right),$
\BE\label{eqW}
W_g (s,t) = \frac{1}{s} \int_{B(x_0,r)\cap B(t,s)} (f(x)-P_{f,x_0} (x-x_0)) \psi\left( \frac{x-t}{s} \right)dx. \nonumber
\EE
If $x\in B(t,s)$ then $|x-x_0| \leq |x-t|+|t-x_0| \leq s + \frac{r}{2}<r$. Hence $B(t,s)\subseteq B(x_0,r)$ and
$$
W_g (s,t) = \frac{1}{s} \int_{B(t,s)} f(x) \psi\left(\frac{x-t}{s}\right)dx - \frac{1}{s} \int_{B(t,s)} P_{f,x_0} (x-x_0) \psi\left(\frac{x-t}{s}\right)dx.$$
The function $\psi$ has $r_{\psi}+1\geq \alpha+1$ vanishing moments and the polynomial $P_{f,x_0}$ is of degree less than $\alpha$, hence
$$\frac{1}{s} \int_{B(t,s)} P_{f,x_0} (x-x_0) \psi\left(\frac{x-t}{s}\right)dx =0.
$$
Therefore,
\BE\label{eqtfi}
\forall s \in \left[ 0, \frac{r}{2} \right], \quad \forall t \in B \left( x_0, \frac{r}{2} \right), \qquad  W_g (s,t) =  \frac{1}{s} \int_{\RR} f(x) \psi\left( \frac{x-t}{s} \right)dx = W_f (s,t).
\EE
Applying \eqref{eqtfi} with $r= 2a$ for $a$ small enough, and using \eqref{caracLp}, we found an upper bound for the $L^p$-leaders \eqref{defplead} by writing
\BENE
\lVert g \rVert_{L^p} & \geq & C_2 \left( \int_{B(x_0,a)} \left(\int_{0}^{a} |W_g (s,t)|^2 \frac{ds}{s} \right)^{\frac{p}{2}} dt \right)^{\frac{1}{p}} \\
& = & C_2 \left( \int_{B(x_0,a)} \left(\int_{0}^{a} |W_f (s,t)|^2 \frac{ds}{s} \right)^{\frac{p}{2}} dt \right)^{\frac{1}{p}} \\
& = & C_2 L_f^{(p)} (a,x_0) a^{\frac{1}{p}}.
\EENE
By applying \eqref{eqnorm}, we deduce that there exists $C'>0$ such that
$$L_f^{(p)} (a,x_0) \leq\frac{1}{C_2} \lVert g\rVert_{L^2}a^{-\frac{1}{p}} \leq C'r^{\alpha + \frac{1}{p}} a^{-\frac{1}{p}} \leq C'a^{\alpha}.$$
\end{proof}

\subsection{Proof of point 2. of the Theorem \ref{propcar}}\label{subsec4.2}

In this section, we proof if $\alpha\notin\NN$ and \eqref{hypL} is verified, then $f$ belongs to $\TPA$.

\begin{proof}
Suppose that for $\alpha\notin\NN$, there exists $C_f>0$ such that for all $a\in\RR_+^*$,
$$|L_f^{(p)} (a,x_0)|\leq C_f a^\alpha.$$
Let us show that $f\in \TPA$. 

Since $\psi$ is even, the admissibility condition \eqref{AdFour} is verified. By the hypothesis (2.4.6) of \cite{Daubechies92} or (B.26) of \cite{Jaffard96b}, we fix $\phi$ a wavelet different from $\psi$ of class $C^{N_\psi} (\RR)$ where $N_\psi>\max(\alpha, 0)$ and compactly supported on $[-1,1]$ with at least one vanishing moment and thus \eqref{eqrerecomp} is valid with this new wavelet.

We separate the cases $ \alpha >0$ and $\alpha \in (-1/p,0]$. 

\medskip

$\bullet$ \textbf{For $\alpha>0$} : Let $r\in (0,1)$ and $x\in B(x_0,r)$. Let
$$
    \CCC (x_0,r) = \{(s,t)\in\RR_+^*\times\RR : |t-x_0|<s+2r\}
$$
For all $(s,t)\in\RR_+^*\times\RR$, we define
\BENE
V_1 (s,t) &=& \left\{\begin{tabular}{ll}
     $W_f (s,t)$ & if $(s,t)\in D_1(x_0,r)=\CCC(x_0,r)\cap \left((0,r)\times \RR\right)$  \\
     $0$ & else
\end{tabular}\right. \\
V_2 (s,t) &=& \left\{\begin{tabular}{ll}
     $W_f (s,t)$ & if $(s,t)\in D_2(x_0,r)=\CCC(x_0,r)\cap \left([1,r)\times \RR\right)$  \\
     $0$ & else
\end{tabular}\right. \\
V_3 (s,t) &=& \left\{\begin{tabular}{ll}
     $W_f (s,t)$ & if $(s,t)\in D_3(x_0,r)=\CCC(x_0,r)\cap \left([r,+\infty)\times \RR\right)$  \\
     $0$ & else
\end{tabular}\right. \\
V_4 (s,t) &=& \left\{\begin{tabular}{ll}
     $W_f (s,t)$ & if $(s,t)\in D_4(x_0,r)= \left(\RR_+^*\times\RR\right)\backslash\CCC(x_0,r)$ \\
     $0$ & else.
\end{tabular}\right.
\EENE

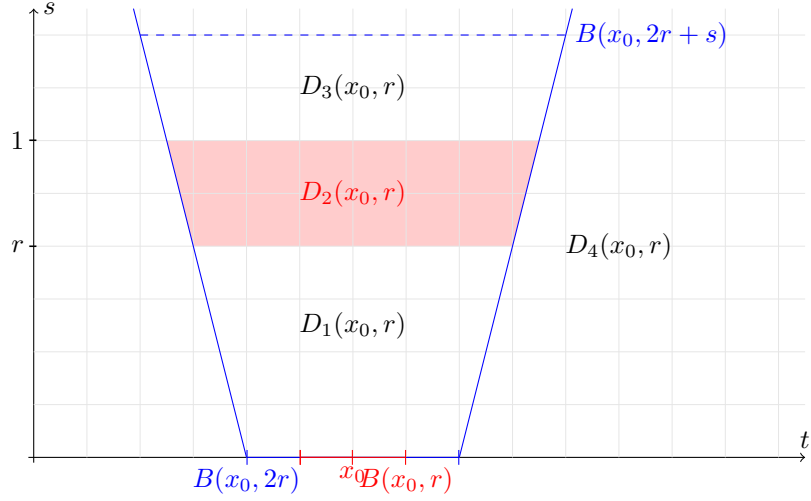
\begin{figure}
    \centering
    \begin{tikzpicture}[scale=0.7]
    \fill[red!20] (2.5,6) -- (9.5,6) -- (9,4) -- (3,4) -- cycle;
    \draw[->,black] (-0.1,0) -- (14.5,0)  node[above] {$t$};
    \draw[->,black] (0,-0.1) -- (0,8.5)  node[right] {$s$};
    \draw[-,very thin,gray!20] (1,0) -- (1,8.5);
    \draw[-,very thin,gray!20] (2,0) -- (2,8.5);
    \draw[-,very thin,gray!20] (3,0) -- (3,8.5);
    \draw[-,very thin,gray!20] (4,0) -- (4,8.5);
    \draw[-,very thin,gray!20] (5,0) -- (5,8.5);
    \draw[-,very thin,gray!20] (6,0) -- (6,8.5);
    \draw[-,very thin,gray!20] (7,0) -- (7,8.5);
    \draw[-,very thin,gray!20] (8,0) -- (8,8.5);
    \draw[-,very thin,gray!20] (9,0) -- (9,8.5);
    \draw[-,very thin,gray!20] (10,0) -- (10,8.5);
    \draw[-,very thin,gray!20] (11,0) -- (11,8.5);
    \draw[-,very thin,gray!20] (12,0) -- (12,8.5);
    \draw[-,very thin,gray!20] (13,0) -- (13,8.5);
    \draw[-,very thin,gray!20] (14,0) -- (14,8.5);
    \draw[-,very thin,gray!20] (0,1) -- (14.5,1);
    \draw[-,very thin,gray!20] (0,2) -- (14.5,2);
    \draw[-,very thin,gray!20] (0,3) -- (14.5,3);
    \draw[-,very thin,gray!20] (0,4) -- (14.5,4);
    \draw[-,very thin,gray!20] (0,5) -- (14.5,5);
    \draw[-,very thin,gray!20] (0,6) -- (14.5,6);
    \draw[-,very thin,gray!20] (0,7) -- (14.5,7);
    \draw[-,very thin,gray!20] (0,8) -- (14.5,8);

    \draw[-,blue] (4,0) -- (1.875,8.5);
    \draw[-,blue] (8,0) -- (10.125,8.5);
    \draw[dashed,blue] (2,8) -- (10,8) node[right] {$B(x_0,2r+s)$};

    \draw[|-|,red] (5,0) -- (6,0)  node[below] {$x_0$};
    \draw[|-|,blue] (8,0) -- (4,0) node[below] {$B(x_0,2r)$};
    \draw[|-|,red] (5,0) -- (7,0) node[below] {$B(x_0,r)$};

    \draw (0,4) node[left] {$r$};
    \draw (0,6) node[left] {$1$};
    \node[] at (0,4) {-};
    \node[] at (0,6) {-};

    \draw[] (6,2.5) node[] {$D_1 (x_0,r)$};
    \draw[red] (6,5) node[] {$D_2 (x_0,r)$};
    \draw[] (6,7) node[] {$D_3 (x_0,r)$};
    \draw[] (11,4) node[] {$D_4 (x_0,r)$};
    
    \end{tikzpicture}
    \caption{Cutting into 4 parts}
    \label{fig:schema}
\end{figure}

We decompose $W_f$ into the four integrals $W_f=V_1+V_2+V_3+V_4$ (see Figure \ref{fig:schema}). For $j=1,2,3,4$, let
$$
    I_j (x) = \frac{1}{c_\phi} \int_{s=0}^{+\infty} \int_{t=-\infty}^{+\infty} V_j (s,t) \left[ \phi\left( \frac{x-t}{s}\right) - \sum_{k=0}^{\lfloor \alpha \rfloor} \frac{\phi^{(k)} \left(\frac{x_0 - t}{s}\right) }{k!} \left( \frac{x-x_0}{s}\right)^k \right] \frac{dtds}{s^2}.
$$
Thus,
$$
I_j (x) = \frac{1}{c_\phi} \int_{D_j (x_0,r)} W_f (s,t) \left(\phi\left(\frac{x-t}{s}\right)-\sum_{k=0}^{\lfloor \alpha \rfloor} \frac{\phi^{(k)} \left(\frac{x_0 - t}{s}\right) }{k!} \left( \frac{x-x_0}{s}\right)^k \right) \frac{dtds}{s^2} 
$$
We will bound the quantities, $\lVert I_j \mathds{1}_{B(x_0,r)}\rVert_{L^p}$ for $j=1,2,3,4$ in order to show that the Taylor-Lagrange polynomial does exist.

\medskip

\textbf{Estimation of $I_1$ :} 
The proof of \eqref{caracLp} provides the upper-bound of the norm $L^p$ of $I_1$ only based on the reconstruction formula \eqref{eqrerecomp} (indeed, this demonstration does not use that $W_f$ is a continuous wavelet transform). We deduce that
$$\lVert I_1 \mathds{1}_{B(x_0,r)} \rVert_{L^p}  \leq  C_1 \left[ \int_0^{+\infty} \left( \int_0^{+\infty} |V_1(s,t)|^2 \frac{ds}{s}\right)^{\frac{p}{2}} dt \right]^{\frac{1}{p}}$$
Since $0<s<r$, one has $B(x_0,2r+s) \subseteq B(x_0,3r)$ and so
\BENE
\lVert I_1 \mathds{1}_{B(x_0,r)} \rVert_{L^p} &\leq &  C_1 \left[ \int_{B(x_0, 3r)} \left( \int_0^{r} |W_{f} (s,t)|^2 \frac{ds}{s}\right)^{\frac{p}{2}} dt \right]^{\frac{1}{p}}\\
& \leq & C_1 (3r)^{\frac{1}{p}} L_f^{(p)} (3r,x_0).
\EENE
From \eqref{hypL}, we conclude that there exists $C_1'>0$ such that
\BE\label{eq1}\lVert I_1 \mathds{1}_{B(x_0,r)} \rVert_{L^p} = C_1 (3r)^{\frac{1}{p}} (3r)^{-\alpha} \leq C'_1 r^{\alpha+\frac{1}{p}}.\EE

\textbf{Estimation of $I_2$ :} According to the Taylor-Lagrange theorem applied to $\phi\in C^{\lfloor\alpha \rfloor + 1}$, there exists $C>0$ such that
$$|I_2 (x)| \leq C \int_{s=r}^1 \int_{t\in B(x_0,2r+s)} |W_f (s,t)| \cdot \left|\frac{x-x_0}{s} \right|^{\lfloor\alpha\rfloor+1} \frac{dtds}{s^2}.$$
There exists a unique $J\in\NN^*$ such that $2^{-J}< r \leq 2^{-(J-1)}$. Thus,
$$|I_2 (x)| \leq C |x-x_0|^{\lfloor\alpha\rfloor+1} \sum_{j=1}^{J}  \int_{B(x_0,2r+2^{-(j-1)})} \left(\int_{2^{-j}}^{2^{-(j-1)}} \frac{|W_f (s,t)|}{s^{\frac{1}{2}}} \frac{ds}{s^{\lfloor\alpha\rfloor+\frac{5}{2}}}\right) dt.$$
By Cauchy-Schwarz,
\BENE
|I_2(x)| &\leq & C |x-x_0|^{\lfloor\alpha\rfloor+1} \sum_{j=1}^{J} \int_{B(x_0,2r+2^{-(j-1)})} \left(\int_{2^{-j}}^{2^{-(j-1)}} |W_f (s,t)|^2 \frac{ds}{s}\right)^{\frac{1}{2}}\\
& & \times \left(\int_{2^{-j}}^{2^{-(j-1)}} \frac{ds}{s^{2 \lfloor\alpha\rfloor+5}}\right)^{\frac{1}{2}} dt,
\EENE
so that
\BENE
|I_2 (x)| &\leq & C |x-x_0|^{\lfloor\alpha\rfloor+1} \sum_{j=1}^{J} 2^{(\lfloor\alpha\rfloor + 2)j} \int_{B(x_0,6\cdot 2^{-j})} \left( \int_{0}^{6\times 2^{-j}} |W_f (s,t)|^2 \frac{ds}{s}\right)^{\frac{1}{2}} dt.
\EENE
According to Hölder's inequality, for $q>0$ such that $1/p+1/q=1$,
\BENE
|I_2(x)| &\leq & C |x-x_0|^{\lfloor\alpha\rfloor+1} \sum_{j=1}^{J} 2^{(\lfloor\alpha\rfloor+2)j} \left[ \int_{B(x_0,6\cdot 2^{-j})} \left( \int_{0}^{6\cdot 2^{-j}} |W_f (s,t)|^2 \frac{ds}{s}\right)^{\frac{p}{2}} dt \right]^{\frac{1}{p}} \\
& &\times \left( \int_{B(x_0,6\cdot 2^{-j})} 1^q dt \right)^{\frac{1}{q}} \\
& \leq & C |x-x_0|^{\lfloor\alpha\rfloor+1} \sum_{j=1}^{J} 2^{(\lfloor\alpha\rfloor+2)j}\cdot (6\cdot 2^{-j})^{\frac{1}{p}} \cdot L_f^{(p)}(6\cdot 2^{-j},x_0) \cdot (2\cdot6 \cdot 2^{-j})^{\frac{1}{q}}.
\EENE
From \eqref{hypL}, there exists $C>0$ such that
$$
|L^{(p)}_f (6\cdot 2^{-j},x_0)| = \left[ \frac{1}{6\cdot 2^{-j}} \int_{B(x_0,6\cdot 2^{-j})} \left( \int_{0}^{6\cdot 2^{-j}} |W_f (s,t)|^2 \frac{ds}{s}\right)^{\frac{p}{2}}dt\right]^{\frac{1}{p}} \leq  C (6\cdot 2^{-j})^{\alpha}.
$$
We deduce that
$$|I_2 (x)| \leq C |x-x_0|^{\lfloor\alpha\rfloor+1} \sum_{j=1}^{J} 2^{(\lfloor\alpha\rfloor+2)j} C (6\cdot 2^{-j})^{\alpha+\frac{1}{p}}(2\cdot 6\cdot 2^{-j})^{\frac{1}{q}}.$$
Since $1/p+1/q=1$, there exists $C'>0$ independent of $J$ and $x$ such that for $\alpha\in (0,1)$ 
$$
|I_2 (x)|\leq C' |x-x_0|^{\lfloor\alpha\rfloor+1}\sum_{j=1}^{J} 2^{(1-\alpha+\lfloor\alpha\rfloor)j}.
$$
We conclude that there exists $C_\alpha >0$ such that $|I_2 (x)|\leq C_\alpha |x-x_0|^{\lfloor\alpha\rfloor+1} 2^{(1-\alpha+\lfloor\alpha\rfloor)J}$. But since $x\in B(x_0,r)$ and $r\in (2^{-J},2^{-(J-1)}]$, we have $|I_2(x)| \leq 2 C_\alpha 2^{-\alpha J} \leq 2 C_\alpha r^{\alpha}$.
So there is a constant $C_2=2^{1+1/p} C_\alpha>0$ independent of $r$ such that
\BE\label{eq2}
\lVert I_2 \mathds{1}_{B(x_0,r)}\rVert_{L^p} = \left(\int_{x_0-r}^{x_0+r} |I_2 (x)|^p dx \right)^{\frac{1}{p}} \leq  2C_\alpha r^\alpha (2r)^{\frac{1}{p}} = C_2 r^{\alpha+\frac{1}{p}}.
\EE

\medskip

\textbf{Estimation of $I_3$ :} Since $\phi\in C^\alpha(\RR)$ and $s>1$, there exists $C>0$ such that
$$\left| \phi\left( \frac{x-t}{s} \right) - \sum_{k=0}^{\lfloor \alpha \rfloor} \frac{\phi^{(k)} \left(\frac{x_0 - t}{s}\right) }{k!} \left( \frac{x-x_0}{s}\right)^k \right| \leq C \left|\frac{x-x_0}{s}\right|^\alpha\leq C |x-x_0|^\alpha.$$
We note
$$C'=\frac{C}{c_\phi} \left(\int_{s=1}^{+\infty} \int_{t\in B(x_0,3+s)} |W_f (s,t)|\frac{dtds}{s^2}\right).$$
Thus, for all $x\in B(x_0,r)$, we have
\BENE
|I_3 (x)| &\leq& \frac{C}{c_\phi} \int_{D_j (x_0,r)} |W_f (s,t)| |x-x_0|^\alpha \frac{dtds}{s^2} \\
&\leq& \frac{C}{c_\phi} r^\alpha \left(\int_{s=1}^{+\infty} \int_{t\in B(x_0,2r+s)} |W_f (s,t)|\frac{dtds}{s^2}\right) \\
&\leq & C' r^\alpha
\EENE
from which, as for $I_2$, we deduce directly that there exists $C_3>0$ such that
\BE\label{eq3}
\lVert I_3 \mathds{1}_{B(x_0,r)}\rVert_{L^p} \leq C_3 r^{\alpha+\frac{1}{p}}.
\EE

\textbf{Estimation of $I_4$ :} For all $s\in\RR_+^*$, one knows that $supp(\phi) \subseteq \overline{B(0,1)}$ and so for $(s,t)\mapsto \phi\left((x-t)/s\right)$, $supp(\phi_{s,t}) \subseteq \overline{B(t,s)}$. Thus, if $t\in \RR\backslash B(x_0,2r+s)$, then $t\notin B(x_0,s)$. Moreover, since $x\in B(x_0,r)$, it follows that
$$2r+s\leq |t-x_0| \leq |t-x|+|x-x_0|\leq |t-x|+r$$
and so $|t-x|\geq r+s \geq s$. We conclude that if $t\in \RR\backslash B(x_0,2r+s)$, then $t\notin B(x_0,s)\cup B(x,s)$, and $\phi\left(\frac{x_0-t}{s}\right)=\phi\left(\frac{x-t}{s}\right)=0$. So
\BE\label{eq4}
I_4(x)=0.
\EE

In conclusion, for $i=1,2,3,4$, the functions $x\mapsto I_i (x) \mathds{1}_{B(x_0, r)}$ belong to $\TPA$ and so according to the wavelet reconstruction formula \eqref{eqrerecomp}, we have $I_1+I_2+I_3 = f(x)-P_{f,x_0} (x-x_0)$ with
$$
    P_{f,x_0} (x-x_0) = \sum_{k=0}^{\lfloor \alpha \rfloor} \int_{s=0}^{+\infty}\int_{t=-\infty}^{+\infty} W_f (s,t)  \frac{\phi^{(k)} \left(\frac{x_0 - t}{s}\right) }{k!} \frac{dtds}{s^{k+2}} (x-x_0)^k
$$
Combining \eqref{eq1}, \eqref{eq2}, \eqref{eq3} and \eqref{eq4}, there exists $C>0$ such that
$$\left(\frac{1}{r} \int_{B(x_0,r)} |f(x)-P_{f,x_0} (x-x_0)|^p dx \right)^{\frac{1}{p}} \leq \frac{1}{r^{\frac{1}{p}}} \sum_{k=1}^{4} \lVert I_k \mathds{1}_{B(x_0,r)}\rVert_{L^p} \leq  C r^{\alpha}.$$
This shows that $f\in\TPA$, i.e. item 1. of Theorem \ref{propcar} for $\alpha>0$.

\medskip

$\bullet$ \textbf{For $\alpha\in (-1/p,0]$} : Let $r\in (0,1)$ and $x\in B(x_0,r)$. Since $\alpha \leq 0$, we can take $P_{f,x_0}=0$, hence
$$f(x) =\frac{1}{c_{\phi}} \int_{s=0}^{+\infty} \int_{t=-\infty}^{+\infty} W_f (s,t) \phi\left(\frac{x-t}{s}\right)\frac{dtds}{s^2}.$$
We have $f(x) = I'_1 (x) + I'_2 (x) + I'_3 (x) + I'_4 (x)$, with for $j=1,2,3,4$, 
$$
I'_j (x)  = \frac{1}{c_{\phi}}  \int_{D_j (x_0,r)} W_f (s,t) \phi\left(\frac{x-t}{s}\right)\frac{dtds}{s^2}.
$$

We need to find an upper-bound for the $L^p$-norm of $x\mapsto |f(x)|\mathds{1}_{B(x_0,r)}(x)$.

\medskip

\textbf{Estimation of $I'_1$ :} The estimate leading to \eqref{eq1} remains valid.

\medskip

\textbf{Estimation of $I'_2$ :} Bounding $\phi$ by $\lVert\phi\rVert_\infty$ since it is continuous compactly supported, we obtain
$$|I'_2 (x)| \leq \frac{\lVert\phi\rVert_\infty}{c_{\phi}} \int_{s=r}^1 \int_{t\in B(x_0,2r+s)} |W_f (s,t)| \frac{dtds}{s^2}.$$
Let $J\in\NN$ such that $2^{-J}< r \leq 2^{-(J-1)}$. Thus, there exists $C_2>0$ such that
$$|I'_2 (x)| \leq C_2 \sum_{j=1}^{J}  \int_{B(x_0,2r+2^{-(j-1)})} \left(\int_{2^{-j}}^{2^{-(j-1)}} \frac{|W_f (s,t)|}{s^{\frac{1}{2}}} \frac{ds}{s^{\frac{3}{2}}}\right) dt.$$
By Cauchy-Schwarz inequality, 
$$|I'_2(x)| \leq C_2  \sum_{j=1}^{J} \int_{B(x_0,2r+2^{-(j-1)})} \left(\int_{2^{-j}}^{2^{-(j-1)}} |W_f (s,t)|^2 \frac{ds}{s}\right)^{\frac{1}{2}} \left(\int_{2^{-j}}^{2^{-(j-1)}} \frac{ds}{s^3}\right)^{\frac{1}{2}} dt.$$
The remaining part of the computation is the same as in the case $\alpha\in (0,1)$ but with $\alpha\in (-1/p,0)$, so there exists $C_4>0$ such that
$$|I'_2 (x)| \leq C_4 \sum_{j=1}^{J} 2^{-\alpha j}.$$
Since $\alpha<0$, there exists $C_\alpha >0$ such that $ |I'_2 (x)| \leq C_\alpha 2^{-\alpha J}$.
Therefore there exists a constant $C>0$ independent of $r$ such that
\BE
\lVert I'_2 \mathds{1}_{B(x_0,r)}\rVert_{L^p}= \left(\int_{B(x_0,r)} |I'_2(x)|^pdx\right)^{\frac{1}{p}} \leq C r^{\alpha+\frac{1}{p}}. \nonumber
\EE

\textbf{Estimation of $I'_3$ :} The reconstruction formula by the continuous wavelet transform restricted to $s\geq 1$ gives a function $\widetilde{g}$ defined by
$$\widetilde{g} (x) = \frac{1}{c_{\phi}} \int_{s=1}^{+\infty}\int_{t\in\RR} W_f (s,t)\phi\left(\frac{x-t}{s}\right) \frac{dtds}{s^2}.$$
The function $\widetilde{g}$ is of class $C^N$ where $N$ is the regularity of the wavelet $\phi$ \cite{Daubechies92}. Therefore, since $\alpha\in (-1/p,0)$, $1\leq r^\alpha$ and therefore there exists $C'>0$ such that
$$C'=\frac{\lVert\phi\rVert_\infty}{c_{\phi}}  \left(\int_{s=1}^{+\infty} \int_{t\in\RR} |W_f (s,t)|\frac{dtds}{s^2}\right)<+\infty.$$
Thus, for all $x\in B(x_0,r)$, we have
$$|I'_3 (x)| \leq \frac{\lVert\phi\rVert_\infty}{c_{\phi}}  \left(\int_{s=1}^{+\infty} \int_{t\in\RR} |W_f (s,t)|\frac{dtds}{s^2}\right) \leq C' \leq C' r^\alpha. $$
And it follows that there exists $C_3>0$ such that
\BE
\lVert I'_3 \mathds{1}_{B(x_0,r)}\rVert_{L^p} \leq C_3 r^{\alpha+\frac{1}{p}}. \nonumber
\EE

\medskip

\textbf{Estimation of $I'_4$ :} It is identical to the $\alpha>0$ case.

\medskip

Combining the 4 previous results, gives $C>0$ such that
$$
\left(\frac{1}{r} \int_{B(x_0,r)} |f(x)-P_{f,x_0}(x-x_0)|^p dx \right)^{\frac{1}{p}} \leq \frac{1}{r^{\frac{1}{p}}} \sum_{k=1}^{4} \lVert I'_k \mathds{1}_{B(x_0,r)}\rVert_{L^p} \leq  C r^{\alpha}.
$$
We conclude that $f\in \TPA$.
\newline
\end{proof}

\section{\texorpdfstring{$p$}{p}-spectrum of sums of random pulses}\label{sec5}

Definition \ref{defpulse} yields locally bounded sample paths of sums of random pulses for $\alpha>0$ and $\eta\in (0,1/\alpha)$ \cite{Saes20}. However, in this part, it is also relevant to consider the non-locally bounded case with $\alpha<0$. We will start by giving sufficient conditions for convergence of such series in $L^p$. More precisely, we will show in a first step that, when $\alpha<0$, $\eta\in (0,1)$ and $\eta-1<\alpha\eta$, the sums of random pulses are in a space $L^p_{loc}(\RR)$ for $p\in [1,-1/(\alpha\eta)+1/\alpha)$. In a second step, we will compute their $p$-multifractal spectra. 

\subsection{Recalls}\label{subsec5.1}

We recall the notations \cite{Saes20}. Let
\BENE
    A_j &=& \{ n\in\NN^* : 2^{j-1} \leq B_n^{\frac{1}{\eta}} < 2^j\} \quad \mbox{if } j\neq 0 \\
    A_0 &=& \{n\in\NN^* : 0 \leq B_n^{\frac{1}{\eta}} < 1\}
\EENE
We also state preliminary lemmas whose proofs can be found in \cite{Saes20} (see Lemmas 3.1, 3.2 and 3.3) see also \cite{Saes21}.

\BL\label{lem1}
Almost surely for $j$ large enough and $\ep_j = \log_2 (j)/(\eta j)$.
$$2^{\eta j (1+\ep_j)} \leq \Card (A_j) \leq 2^{\eta j (1+\ep_j)}.$$
\EL

\BL\label{lem2}
Almost surely, there exists $K_1,K_2>0$ such that for all $j\in\NN^*$ and $n\in A_j$,
$$K_1 2^{\eta j (1-\ep_j)} \leq B_n, C_n \leq K_2 2^{\eta j (1+\ep_j)}.$$
\EL

For all $x,r\in [0,1]$, we note $T_n (x,r)=1$ if $B(X_n,B_n^{1/\eta})\cap B(x,r)\neq \emptyset$ and $T_n (x,r)=0$ else. For $r=0$, we write $T_n(x)=T_n(x,0)$.

\BL\label{lem3}
Almost surely, there exists $K>0$ such that for all $x\in [0,1]$ and $J,j\in\NN$ with $j\neq 0$,
$$\sum_{n\in A_j} T_n (x,2^{-\eta J}) \leq K j^2 \max\{1, 2^{\eta (j-J)}\} \quad \mbox{et} \quad \sum_{n\in A_j} T_n (x) \leq Kj^2.$$
\EL

Let $p_0$ be a sufficiently large integer such that $\rho > \frac{3-3\alpha}{1-\alpha\eta}$. Let $gamma \in (0,\frac{1}{\eta}-1)$ and consider for all $j\in\NN$, the sets
$$\widetilde{A_j} = \bigcup_{j'=\lfloor (1-\rho \eta\ep_j)j \rfloor}^{\lfloor \gamma j \rfloor} A_j$$
and 
\BE\label{defI}
\mathcal{I}_j = \{n\in A_j : \forall m \in\widetilde{A_j},\quad n\neq m,\quad B(X_n,B_n^{-\frac{1}{\eta}})\cap B(X_m,B_m^{-\frac{1}{\eta}}) = \emptyset\}. \nonumber
\EE
We define the families of sets $G_\delta$ and $G'_{\delta}$ by
\BE
G_\delta  = \limsup_{j\rightarrow +\infty} \bigcup_{n\in A_j} B(X_n,B_n^{-\delta }) \nonumber
\EE
and
\BE
G'_\delta  = \limsup_{j\rightarrow +\infty} \bigcup_{n\in \widetilde{A}_j } B(X_n,B_n^{-\delta(1-\widetilde{\ep}_j)}) \quad \mbox{and} \quad  \nonumber 
\EE

We prove the Theorem \ref{pthm1} which yields the multifractal analysis of random sums of pulses in the $p$-exponent framework.

When $p=+\infty$, this result boils down to the results proved in \cite{Saes20}. Let us notice that the sample paths of $F_{\alpha,\eta}$ have the same $p$-spectrum as the random lacunar wavelets series studied \cite{Leonarduzzi16}.

\subsection{Suitable of \texorpdfstring{$p$}{p} values for \texorpdfstring{$p$}{p}-multifractal analysis}\label{subsec5.2}

First, let us determine the spaces $L^p_{loc}(\RR)$ to which the function $f$ belongs according to the values of $\alpha<0$.

\BP\label{propFLp}
Let $\psi:\RR\rightarrow\RR$ a continuous lipschitzian function with support on $[-1,1]$. Let $\eta\in (0,1)$, $\alpha <0$ such that $\eta -1 < \alpha \eta$. Almost surely, for all $p\in (1,-1/(\alpha\eta) + 1/\alpha)$, the sample paths of the process $F_{\alpha,\eta}:\RR\rightarrow\RR$ defined by \eqref{defP} belongs to $ L^{p}_{loc} ([0,1])$.
\EP
Note that the condition $\eta -1 < \alpha \eta$ guarantees that $-1/(\alpha\eta) + 1/\alpha >1$.

\begin{proof}
Let $p\in [1,-1/(\alpha\eta) + 1/\alpha)$, and consider
\BE\label{eqpr1}
F_{\alpha,\eta}(x)=\sum_{j=0}^{+\infty}  F_j (x) \ \mbox{ where } \ F_j (x)=\sum_{n\in A_j}C_n^{-\alpha} \psi (B_n^{\frac{1}{\eta}} (x-X_n)). \nonumber
\EE
We estimate separately the $L^p$-norm of each $F_j$. Applying the definition of $T_n$ and the Lemma \ref{lem3}, there exists $K>0$ such that for all $x\in\RR$,
$$|F_j (x)|\leq K \lVert\psi\rVert_\infty \max_{n\in A_j} \{C_n^{-\alpha}\} \sum_{n\in A_j} T_n(x) \leq K j^2 \max_{n\in A_j} \{C_n^{-\alpha}\} \mathds{1}_{\bigcup_{n\in A_j} B(X_n,B_n^{-1/\eta})} (x)$$
and 
$$\lVert F_j \rVert_{L^p} \leq K j^{2} \max_{n\in A_j} \{C_n^{-\alpha}\} \lVert \mathds{1}_{\bigcup_{n\in A_j} B(X_n,B_n^{-1/\eta})} \rVert_{L^p}.$$
But the Lebesgue measure of the support of $\bigcup_{n\in A_j}\mathds{1}_{ B(X_n,B_n^{-1/\eta})}$ is bounded by $\Card(A_j) \max_{n\in A_j}$ $\{B_n^{-1/\eta}\}$, we deduce that
$$\lVert F_j\rVert_{L^p} \leq Kj^2 \max_{n\in A_j}\{C_n^{-\alpha}\} \Card(A_j)^{1/p}\max_{n\in A_j}\{B_n^{-\frac{1}{\eta p}}\}.$$
By Lemmas \ref{lem1} and \ref{lem2},
$$\lVert F_j \rVert_{L^p} \leq Kj^2 j^{\alpha}2^{-\alpha\eta j} j^{1/p} 2^{\frac{\eta}{p}j}  j^{1/\eta p} 2^{-\frac{1}{p}j} \leq K j^{C(\alpha, \eta, p)} 2^{(-\alpha\eta +\frac{\eta}{p} - \frac{1}{p})j}.$$
Finally,
$$\sum_{j\in\NN} \lVert F_j \rVert_{L^p} \leq K\sum_{j\in\NN} j^{C(\alpha, \eta, p)} 2^{(-\alpha\eta +\frac{\eta}{p} - \frac{1}{p})j}.$$
By hypothesis, $-\alpha\eta+\eta/p +1/p <0$. We conclude that almost surely the series $\sum\limits_{j\in\NN} |F_j(x)|$ converges in $L^p$ for $p\in [1,-1/(\alpha\eta) + 1/\alpha)$. So $\sum F_j$ converges too and it converges to $F_{\alpha,\eta}$. Thus, $F_{\alpha,\eta}\in L^p$ for $p\in [1,-1/(\alpha\eta) + 1/\alpha)$, as stated.
\end{proof}

\subsection{Study of \texorpdfstring{$p$}{p}-exponents}\label{subsec5.3}

\BP\label{pprop1}
Let $\psi:\RR\rightarrow\RR$ be a continuous Lipschitz function with support on $[-1,1]$. Let $\eta\in (0,1)$, $\alpha <0$ such that $\eta -1 < \alpha \eta$. Consider the process $F_{\alpha,\eta}$ defined by \eqref{defP}. Almost surely, for all $p\in [1,-1/(\alpha\eta) + 1/\alpha)$, for all $\delta\in (1,\frac{1}{\eta} )$ and for all $x_0\notin G_\delta$,
$$
    h^{(p)}_{F_{\alpha,\eta}} (x_0) \geq \frac{\alpha}{\delta}+\frac{1-\delta\eta}{\delta\eta p}  =\frac{\alpha+1/(\eta p)}{\delta}+1/p.
$$
\EP

\begin{proof}
Let $x\notin G_\delta$. By definition, there exists $J_x\in\NN$ depending on $x$ such that for any integer $j\geq J_x$, for any $n\in A_j$, $x\notin B(X_n,B_n^{-\delta})$. The sum $$\sum_{j=1}^{J_x} F_j$$ has a finite number of terms, so its global regularity is that of the wavelet. So, in estimating the regularity of $F_{\alpha,\eta}$ in $x$, we can assume that the sum in $F_{\alpha,\eta}$ is taken to $j\geq J_x$.

We take $r>0$ fixed. To estimate $\lVert F_{\alpha,\eta} \rVert_{L^p (B(x,r))}$, we decompose $F_{\alpha,\eta}$ into $\sum_{j\in\NN} F_j$ as in \eqref{eqpr1} and we will determine first for which $j$ the supports of the pulses indexed by elements of $A_j$ can intersect $B(x,r)$. Let $j\geq J_x$, $n\in A_j$. Since $x\notin G_\delta$,
$$|x-X_n| \geq B_n^{-\delta} \geq 2^{-\delta\eta j}.$$
If $B(x,r)$ intersects $B(X_n,B_n^{-1/\eta})$, it implies that
$$r\geq \frac{1}{2} B_n^{-\delta} \geq \frac{1}{2} 2^{-\eta\delta j}.$$
Let $J\in\NN$ be the first integer such that $2^{-\eta J} \leq (2r)^{1/\delta}$.

$\bullet$ If $B_n\in A_j$ for $j<J$, then $B(X_n,B_n^{-\delta}) \cap B(x,r) = \emptyset$ and so $\lVert F_j \rVert_{L^p (B(x,r))} =0$.

$\bullet$ If $B_n\in A_j$ for $j\geq J$, then Lemma \ref{lem3} states that there are at most $Kj^2 2^{\eta (j-J)}$ pulses in $A_j$ whose support intersects $B(x,r)$. Moreover, at a given point, at most $j^2$ overlap. A computation similar to the proof of Proposition \ref{propFLp} gives
\[\int_{B(x,r)} |F_j (x)|^p dx \leq C j^{2p} K j^2 2^{\eta (j-J)} \left(\max_{n\in A_j} C_n^{-\alpha}\right)^p \max_{n\in A_j} B_n^{-\frac{1}{\eta}}.\]
By Lemma \ref{lem2}, 
$$\int_{B(x,r)} |F_j (x)|^p dx \leq C j^{2p} Kj^2 2^{\eta (j-J)} 2^{-\eta p \alpha j} 2^{-j} 2^{-\alpha\eta p j \ep_j}.$$
Since $2^{\eta j \ep_j}=j$,
$$\lVert F_j \rVert_{L^p (B(x,r))} \leq C j^{2+\frac{2}{p}-\alpha} 2^{(-\eta\alpha+\frac{\eta-1}{p})j} 2^{-\frac{\eta}{p}J}.$$
The series of norms converges since $-\alpha\eta+\frac{\eta-1}{p}<0 \Longleftrightarrow p < \frac{1}{\alpha}-\frac{1}{\alpha\eta} $ and thus
$$\exists c_1>0, \quad  \sum_{j=J}^{+\infty} \lVert F_j\rVert_{L^p (B(x,r))} \leq J^{c_1} 2^{(-\eta\alpha - \frac{1}{p})J} \leq 2^{-(\alpha\eta + \frac{1}{p})J} 2^{c_1\eta J \ep_J}.$$
Since $2^{-J} \leq C r^{1/(\delta\eta)}$,
$$\sum_{j=J}^{+\infty} \lVert F_j\rVert_{L^p (B(x,r))} \leq C r^{\frac{\alpha}{\delta} + \frac{1}{p\delta\eta}} w(r)$$
with $w(r)=r^{\ep(r)}$ where $\lim_{r\rightarrow 0} \ep (r) =0$.
Paying attention to the fact that the $p$-exponent is in fact given by 
$$h^p_{F_{\alpha,\eta}} (x) = \liminf_{r\rightarrow 0} \frac{\log \left( \frac{1}{r} \int_{B(x,r)} | F_{\alpha,\eta}(t)-P_{x} (t-x) dt |^p \right)^{\frac{1}{p}}}{\log\left(\frac{1}{r}\right)}$$
and that for $\alpha<0$, we have $P_{x}=0$, we obtain by taking into account the factor $1/r$ before the integral 
$$h^p_{F_{\alpha,\eta}} (x)  \geq \frac{\alpha}{\delta} +\frac{1}{p\delta\eta} + \frac{1}{p}= \frac{\alpha+\frac{1}{\eta p}}{\delta}+\frac{1}{p}.$$
\end{proof}

\BP\label{pprop2}
Let $\psi:\RR\rightarrow\RR$ be a continuous Lipschitz function with support on $[-1,1]$. Let $\eta\in (0,1)$, $\alpha <0$ such that $\eta -1 < \alpha \eta$. Let $F_{\alpha,\eta}$ be defined by \eqref{defP}. Almost surely, for all $p\in [1,-1/(\alpha\eta) + 1/\alpha)$, for all $\delta\in (1,\frac{1}{\eta} )$ and for all $x_0\notin G'_\delta$,
$$
    h^{(p)}_{F_{\alpha,\eta}} (x_0) \leq \frac{\alpha}{\delta}+\frac{1-\delta\eta}{\delta\eta p}  =\frac{\alpha+1/(\eta p)}{\delta}+1/p.
$$
\EP

\begin{proof}
We choose the wavelet function $\phi$ to compute the continuous wavelet transform to belong to $C^1 (\RR)$, even and supported in $[-1,1]$ with at least one vanishing moment. We also impose that
$$\int_{-1}^{1} \phi(u)\psi(u)du \neq 0$$
and by continuity of the continuous wavelet transform, there exists $c\in (0,1)$ such that for all $\ep\in [-c,c]$,
\BE\label{propphi}
\int_{-1}^{1} \phi\left( \frac{u-\ep}{1+\ep} \right) \psi(u)du \neq 0.
\EE
The existence of such a $\phi$ function is straightforward.

The continuous wavelet transform of the sums of random pulses $F_{\alpha,\eta}$, is, for any $(s,t)\in\RR_+^* \times \RR$
$$W_{F_{\alpha,\eta}} (s,t)  = \frac{1}{s} \int_{\RR} F_{\alpha,\eta} (x) \phi\left( \frac{x-t}{s}\right)dx = \sum_{n\in\NN^*} C_n^{-\alpha} d_n (s,t),$$
where 
$$d_n (s,t) = \frac{1}{s} \int_{\RR} \psi_n (x) \phi\left(\frac{x-t}{s}\right)dx \quad \mbox{with} \quad \psi_n (x) = \psi ( B_n^{\frac{1}{\eta}} (x-X_n)).$$
The lemma below, proved in \cite{Demichel06} (Proposition 2.2.1), allows to obtain a upper-bound for the wavelet coefficients of $F_{\alpha,\eta}$.  

\BL\label{lem4} If $\phi$ has zero integral, then there exists $K>0$ such that
\BE
\forall (s,t)\in \RR_+^*\times\RR,\ |d_n(s,t)|\leq K \min\{ s B_n^{\frac{1}{\eta}},s^{-1} B_n^{-\frac{1}{\eta}} \} T_n(t,s). \nonumber
\EE
\EL

Let $x_0\in G'_\delta$. There exists sequences of integers $(n_k)_{k\in\NN}$ and $(j_k)_{k\in\NN}$ such that $n_k\in \mathcal{I} _{j_k}$ and $x_0\in\bigcap\limits_{n\in \mathcal{I} _{j_k}} B(X_n,B_n^{-\delta(1-\widetilde{\ep_{j_k}})})$.

Let $k\in\NN^*$ with $n_k\in \mathcal{I}_{j_k}$ and $\ep\in [-c,c]$. We estimate the continuous wavelet transforms $W_{F_{\alpha,\eta}} (B_{n_k}^{-\frac{1}{\eta}}+\ep B_{n_k}^{-\frac{1}{\eta}}, X_{n_k}+\ep B_{n_k}^{-\frac{1}{\eta}})$.

Let $J_k=\lfloor(1-\rho\eta\ep_{j_k})j_k \rfloor$ and $\widetilde{J}_k= \lfloor \gamma j_k \rfloor$, where $\gamma\in (0,1/\eta)$ has been set. We decompose $W_F (B_{n_k}^{-\frac{1}{\eta}}+\ep B_{n_k}^{-\frac{1}{\eta}}, X_{n_k}+\ep B_{n_k}^{-\frac{1}{\eta}})$ in
\BE
W_F (B_{n_k}^{-\frac{1}{\eta}}+\ep B_{n_k}^{-\frac{1}{\eta}}, X_{n_k}+\ep B_{n_k}^{-\frac{1}{\eta}}) = S_1 + S_2 + S_3. \nonumber
\EE 
with
\BENE
S_1 &=& \sum_{j=0}^{J_k-1}\sum_{n\in A_j}C_n^{-\alpha} d_n(B_{n_k}^{-\frac{1}{\eta}} (1+\ep),X_{n_k}+\ep B_{n_k}^{-\frac{1}{\eta}}), \\
S_2 &=& \sum_{j=J_k}^{\widetilde{J}_k}\sum_{n\in A_j}C_n^{-\alpha} d_n(B_{n_k}^{-\frac{1}{\eta}} (1+\ep),X_{n_k}+\ep B_{n_k}^{-\frac{1}{\eta}}), \\
S_3 &=& \sum_{j=\widetilde{J}_k+1}^{+\infty}\sum_{n\in A_j}C_n^{-\alpha} d_n(B_{n_k}^{-\frac{1}{\eta}} (1+\ep),X_{n_k}+\ep B_{n_k}^{-\frac{1}{\eta}}).
\EENE
The idea is find a lower-bound for $W_F$. To do this, we show that the term $S_2$ is much larger than the two other terms $S_1$ and $S_3$ (which correspond to the high and low frequency terms). 

Let us first consider the sum $S_2$. By the definition of $\mathcal{I}_j$ in \eqref{defI}, there exists a unique $\widetilde{n}_k\in \widetilde{A}_{j_k}$ such that $x_0\in B(X_{\widetilde{n}_k},B_{\widetilde{n}_k}^{-\frac{1}{\eta}})$. Thus $\widetilde{n}_k=n_k$ and
$$
    S_2 = C_{n_k}^{-\alpha} d_{n_k} (B_{n_k}^{-\frac{1}{\eta}} (1+\ep), X_{n_k} + \ep B_n^{-\frac{1}{\eta}}).
$$
We have
\BENE
d_{n_k} (B_{n_k}^{-\frac{1}{\eta}} (1+\ep),X_{n_k}+\ep B_{n_k}^{-\frac{1}{\eta}}) & = & B_{n_k}^{\frac{1}{\eta}} \int_\RR \psi (B_{n_k}^{\frac{1}{\eta}} (x-X_{n_k})) \phi\left( \frac{x-X_{n_k}-\ep B_{n_k}^{-\frac{1}{\eta}}}{(1+\ep)B_{n_k}^{-\frac{1}{\eta}}}\right)dx\\
&=& \int_\RR \psi (u) \phi\left( \frac{B_{n_k}^{-\frac{1}{\eta}}u-\ep B_{n_k}^{-\frac{1}{\eta}}}{(1+\ep)B_{n_k}^{-\frac{1}{\eta}}}\right)du \\
&=& \int_{\RR} \psi (u) \phi\left( \frac{u-\ep}{1+\ep}\right)du.
\EENE
The condition \eqref{propphi} implies that for a constant $K_2>0$ (depending on $\psi$ and $\phi$ only), we have according to the Lemma \ref{lem2} with $\alpha<0$,  
$$
    |S_2|\geq K_2  C_{ {n}_k}^{-\alpha} \geq  K _2  2^{-\alpha\eta (1-\ep_{j_k})j_k} \geq  2^{\alpha (1-\ep_{j_k})}K_2  B_{n_k}^{-\alpha (1-\ep_{j_k})}.
$$
So
\BE\label{minS2}
\exists K_2'>0, \quad |S_2|\geq K'_2 B_{n_k}^{-\alpha (1-\ep_{j_k})}.
\EE

 \smallskip

Then,by Lemmas \ref{lem2} and \ref{lem4} with $\alpha<0$, one has
\BENE
|S_1| & \leq &  \sum_{j=0}^{J_k-1} \sum_{n\in A_j} C_n^{-\alpha} |d_n(B_{n_k}^{-\frac{1}{\eta}} (1+\ep),X_{n_k}+\ep B_{n_k}^{-\frac{1}{\eta}} )|\\
& \leq & \sum_{j=0}^{J_k-1} \sum_{n\in A_j} C_n^{-\alpha}  \min\{ (1+\ep) B_{n_k}^{-\frac{1}{\eta}} B_n^{\frac{1}{\eta}}, (1+\ep)^{-1} B_{n_k}^{\frac{1}{\eta}} B_n^{-\frac{1}{\eta}} \} \\
& & \hspace{2cm} \times\ T_n(X_{n_k}+\ep B_{n_k}^{-\frac{1}{\eta}},(1+\ep) B_{n_k}^{-\frac{1}{\eta}}).
\EENE
Since $\ep\in [-c,c]$ with $c\in (0,1)$,
$$T_n (X_{n_k} + \ep B_{n_k}^{-\frac{1}{\eta}} , (1+\ep) B_{n_k}^{-\frac{1}{\eta}} ) \leq T_n(X_{n_k},2^{-j_k+2}).$$
Therefore, 
$$|S_1|  \leq  \sum_{j=0}^{J_k-1}  2^{-\alpha\eta j(1+\ep_j)} \min\{ (1+\ep) B_{n_k}^{-\frac{1}{\eta}} 2^j, (1+\ep)^{-1} B_{n_k}^{\frac{1}{\eta}} 2^{-j+1} \} \sum_{n\in A_j} T_n(X_{n_k},2^{-j_k+2}).$$
According to Lemma \ref{lem3} and since $j<(1-\rho\eta\ep_{j_k})j_k$,
$$
\exists K_1, K_2 >0,\ |S_1| \leq  K (1+\ep)  B_{n_k}^{-\frac{1}{\eta}} \sum_{j=0}^{J_k-1} j^{2-\alpha} 2^{(1-\alpha\eta) j} \leq  K_1 B_{n_k}^{-\frac{1}{\eta}} j_k^{3-\alpha} 2^{(1-\alpha\eta)(1-\rho\eta\ep_{j_k})j_k}.
$$
Therefore $j_k=2^{\eta\ep_{j_k} j_k}$ and that $n_k\in \mathcal{I} _{j_k}$, $2^{j_k}\leq 2 B_{n_k}^{\frac{1}{\eta}}$, 
$$
    \exists K'_1>0, \quad |S_1| \leq  K_1 B_{n_k}^{-\frac{1}{\eta}} B_{n_k}^{(3-\alpha)\ep_{j_k}} B_{n_k}^{(\frac{1}{\eta}-\alpha)(1-\rho\eta\ep_{j_k})} \leq  K'_1 B_{n_k}^{-\alpha-(\rho-3+\alpha-\alpha\eta \rho)\ep_{j_k}}. 
$$
By our choice of the integer $\rho>\frac{3-3\alpha}{1-\alpha\eta}$, $\rho-3+\alpha-\alpha\eta \rho >-2\alpha$. Thus
\BE
\label{majS1}
    |S_1| \leq  K'_1    B_{n_k}^{-\alpha(1+2\ep_{j_k})}  .
\EE

\smallskip

Now we bound $S_3$ in the same way as $S_1$ with Lemmas \ref{lem2} and \ref{lem4}. Thus
\BENE
|S_3| & \leq &  \sum_{j=\widetilde{J}_k+1}^{+\infty} \sum_{n\in A_j} C_n^{-\alpha} |d_n(B_{n_k}^{-\frac{1}{\eta}} (1+\ep),X_{n_k}+\ep B_{n_k}^{-\frac{1}{\eta}})| \\
& \leq & \sum_{j=\widetilde{J}_k+1}^{+\infty} \sum_{n\in A_j} C_n^{-\alpha} \min\{(1+\ep) B_{n_k}^{-\frac{1}{\eta}} B_n^{\frac{1}{\eta}},(1+\ep)^{-1} B_{n_k}^{\frac{1}{\eta}} B_n^{-\frac{1}{\eta}} \} \\
& & \hspace{2cm} \times\ T_n(X_{n_k}+\ep B_{n_k}^{-\frac{1}{\eta}}, (1+\ep) B_{n_k}^{-\frac{1}{\eta}}) \\
& \leq & \sum_{j=\widetilde{J}_k+1}^{+\infty}  2^{-\alpha\eta j(1-\ep_j)} B_{n_k}^{-\frac{1}{2\eta}} \min\{ B_{n_k}^{-\frac{1}{\eta}} 2^j,B_{n_k}^{\frac{1}{\eta}} 2^{-j+1} \} \sum_{n\in A_j} T_n(X_{n_k},2^{-j_k+2}) .  
\EENE

Therefore $j>\widetilde{J}_k=\lfloor \gamma j_k\rfloor $, the above minimum is reached at $B_{n_k}^{\frac{1}{\eta}} 2^{-j+1}$. Then, according to Lemma \ref{lem3}, there exists a constant $K>0$ such that the sum 
$$\sum_{n\in A_j} T_n (X_n,2^{-j_k+2})$$
is bounded by $Kj^2$ when  $j\leq (j_k-2)/\eta$, and by $Kj^2 2^{\eta (j-j_k/\eta)}$ when $j>(j_k-2)/\eta$. Thus by Lemma \ref{lem2}, there exists a constant $K_3>0$ (which can change value at each line but does not depend on $k$ or other parameters) such that
\BENE
|S_3|& \leq & \frac{K_3}{1+\ep} \left( \sum_{j= \lfloor \gamma j_k\rfloor }^{\lfloor j_k/\eta \rfloor} j^{2-\alpha} 2^{-\alpha\eta j} B_{n_k}^{\frac{1}{\eta}}2^{-j} +   \sum_{j=\lfloor j_k/\eta \rfloor+1}^{+\infty} j^{2-\alpha} 2^{-\alpha\eta j} B_{n_k}^{\frac{1}{\eta}} 2^{-j}2^{\eta \frac{j-(j_k-2)}{\eta}}\right)  \\
& \leq & \frac{K_3}{1+\ep} B_{n_k}^{\frac{1}{\eta}} \left( \sum_{j=\lfloor \gamma j_k\rfloor }^{\lfloor j_k/\eta\rfloor} j^{2-\alpha} 2^{-(1+\alpha\eta)j} + 2^{-j_k+2} \sum_{j=\lfloor j_k/\eta \rfloor+1}^{+\infty} j^{2-\alpha} 2^{(\eta-1-\alpha\eta)j} \right). 
\EENE
Since $(1+\alpha\eta) > \eta >0$ by hypothesis, the first sum above is bounded by 
$$\sum_{j= \lfloor \gamma j_k\rfloor }^{\lfloor j_k/\eta\rfloor} j^{2-\alpha} 2^{-(1+\alpha\eta)j}  \leq K_3  {j}_k^{2-\alpha} 2^{-(1+\alpha\eta)\gamma j_k }$$
and the second sum by
$$
2^{-j_k+2} \sum_{j=\lfloor j_k/\eta \rfloor+1}^{+\infty} j^{2+\alpha} 2^{(\eta-1-\alpha\eta)j}   \leq K_3 2^{-j_k}j_k^{2-\alpha} 2^{(\eta-1-\alpha\eta)\frac{j_k-2}{\eta}} \leq K_3 j_k^{2-\alpha}2^{-\frac{j_k}{\eta}(1+\alpha\eta)}.
$$
As $B_{n_k}^{\frac{1}{\eta}} \sim 2^{j_k}$ and $j_k=2^{j_k \eta \ep_{j_k}} \sim  B_{n_k}^{\ep_{j_k}}$ and $\gamma < 1/\eta $, we obtain that
\BENE
|S_3| & \leq & K_3 B_{n_k}^{\frac{1}{\eta}}  {j}_k^{2-\alpha} 2^{-(1+\alpha\eta)\gamma j_k }  +  K_3 B_{n_k}^{\frac{1}{\eta}} j_k^{2-\alpha}2^{-\frac{j_k}{\eta}(1+\alpha\eta)}\\
&\leq  &   2K_3 B_{n_k} ^{ (2-\alpha)\ep_{j_k} + \frac{1}{\eta}- (\frac{1}{\eta}+\alpha)\gamma } \\
& \leq & 2K_3 B_{n_k}^{(2-\alpha)\ep_{j_k} - \frac{\alpha\eta \gamma +\gamma-1}{\eta}}.
\EENE
We notice that $\frac{\alpha\eta \gamma +\gamma-1}{\eta} - (2-\alpha)\ep_{j_k}> \alpha (1+
2\ep_{j_k})$. Indeed,
$$\frac{\gamma(\alpha\eta + 1)-1}{\eta} > \alpha + (2+\alpha)\ep_{j_k}.$$
Thus, 
\BE\label{majS3}
|S_3| \leq K_3 B_{n_k}^{-\alpha (1+2\ep_{j_k})}.
\EE
\smallskip

By comparing \eqref{minS2}, \eqref{majS1} and \eqref{majS3}, we obtain
\BE\label{eq33}
|W_{F_{\alpha,\eta}} (B_{n_k}^{-\frac{1}{\eta}} (1+\ep),X_{n_k}+\ep B_{n_k}^{-\frac{1}{\eta}})|\geq K B_{n_k}^{-\alpha (1-2\ep_{j_k})}. 
\EE

So for all $\ep\in [-c,c]$, the inequality \eqref{eq33} holds. In particular, there exists $K>0$ such that for all $(s,t) \in [(1-c)B_{n_k}^{-\frac{1}{\eta}}, (1+c) B_{n_k}^{-\frac{1}{\eta}}]\times [X_{n_k}-cB_{n_k}^{-\frac{1}{\eta}}, X_{n_k}+cB_{n_k}^{-\frac{1}{\eta}}]$,
\BE\label{concluW}
|W_{F_{\alpha,\eta}} (s,t)| \geq K B_{n_k}^{-\alpha (1-2\ep_{j_k})}.
\EE
In order to conclude the proof, it remains to estimate 
$$L_{F_{\alpha,\eta}}^{(p)} (B_{n_k}^{-\frac{1}{\eta}}, X_{n_k}) = \left[ B_{n_k}^{\frac{1}{\eta}}\int_{X_{n_k}-B_{n_k}}^{X_{n_k}+B_{n_k}^{-\frac{1}{\eta}}} \left( \int_0^{B_{n_k}^{-\frac{1}{\eta}}} |W_{F_{\alpha,\eta}} (s,t)|^2 \frac{ds}{s} \right)^{\frac{p}{2}} dt \right]^{\frac{1}{p}}.$$

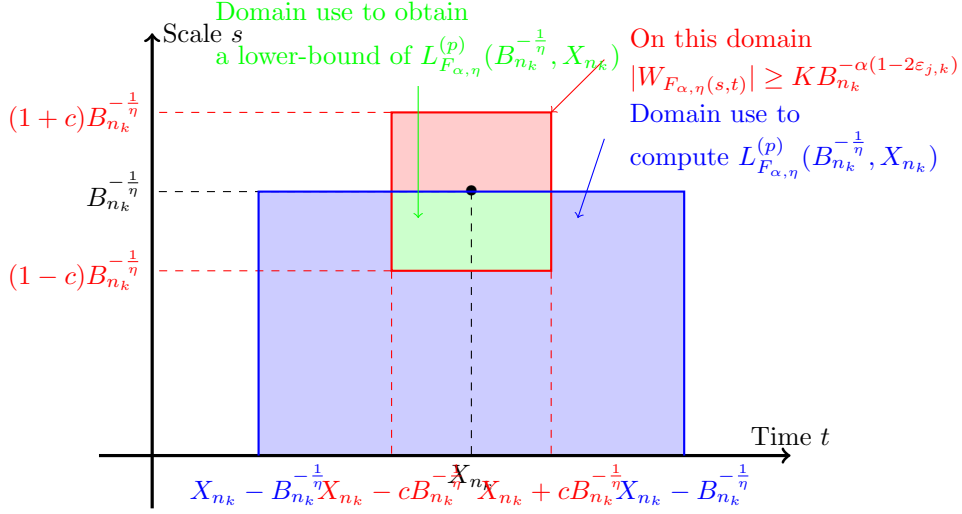
\begin{figure}
\begin{tikzpicture}[scale=0.7]
\fill[blue!20] (2,0) -- (2,5) -- (10,5) -- (10,0) -- cycle ;
\fill[red!20] (4.5,3.5) -- (7.5,3.5) -- (7.5,6.5) -- (4.5,6.5) -- cycle;
\fill[green!20] (4.5,3.5) -- (7.5,3.5) -- (7.5,5) -- (4.5,5) -- cycle;

\draw [->,very thick] (-1,0) -- (12,0) node[above]{Time $t$} ;
\draw [->,very thick] (0,-1) -- (0,8) node[right]{Scale $s$} ;
\draw [dashed] (6,5) -- (6,0) node[below] {$X_{n_k}$};
\draw [dashed] (6,5) -- (0,5) node[left] {$B_{n_k}^{-\frac{1}{\eta}}$};
\draw (6,5) node{$\bullet$};

\draw [dashed,red] (4.5,3.5) -- (4.5,0) node[below] {$X_{n_k}-cB_{n_k}^{-\frac{1}{\eta}}$};
\draw [dashed,red] (7.5,3.5) -- (7.5,0) node[below] {$X_{n_k}+cB_{n_k}^{-\frac{1}{\eta}}$};
\draw [dashed,red] (4.5,6.5) -- (0,6.5) node[left] {$(1+c)B_{n_k}^{-\frac{1}{\eta}}$};
\draw [dashed,red] (4.5,3.5) -- (0,3.5) node[left] {$(1-c)B_{n_k}^{-\frac{1}{\eta}}$};
\draw[thick, red] (4.5,3.5) -- (7.5,3.5) -- (7.5,6.5) -- (4.5,6.5) -- cycle;

\draw[blue] (2,0) node[below]{$X_{n_k}-B_{n_k}^{-\frac{1}{\eta}}$};
\draw[blue] (10,0) node[below]{$X_{n_k}-B_{n_k}^{-\frac{1}{\eta}}$};
\draw[thick, blue] (2,0) -- (2,5) -- (10,5) -- (10,0) ;

\draw[<-,red] (7.5,6.5) -- (8.5,7.5) node[right]{\begin{tabular}{l} On this domain \\ $|W_{F_{\alpha,\eta} (s,t)}| \geq K B_{n_{k}}^{-\alpha(1-2\ep_{j,k})}$\end{tabular}};

\draw[<-,blue] (8,4.5) -- (8.5,6) node[right]{\begin{tabular}{l} Domain use to \\ compute $L_{F_{\alpha,\eta}}^{(p)} (B_{n_k}^{-\frac{1}{\eta}},X_{n_k})$\end{tabular}};

\draw[<-,green] (5,4.5) -- (5,7) node[above]{\begin{tabular}{l} Domain use to obtain \\ a lower-bound of $L_{F_{\alpha,\eta}}^{(p)} (B_{n_k}^{-\frac{1}{\eta}},X_{n_k})$\end{tabular}};
\end{tikzpicture}
\caption{Decomposition time/scale to lower-bound $L_{F_{\alpha,\eta}}^{(p)} (B_{n_k}^{-\frac{1}{\eta}},X_{n_k})$}
\label{FigTimeScale}
\end{figure}

From \eqref{concluW} (which can be visualized using Figure \ref{FigTimeScale}), we have
\BEN
L_{F_{\alpha,\eta}}^{(p)} (B_{n_k}^{-\frac{1}{\eta}}, X_{n_k}) & \geq& \left[ B_{n_k}^{\frac{1}{\eta}}\int_{X_{n_k}-cB_{n_k}}^{X_{n_k}+cB_{n_k}^{-\frac{1}{\eta}}} \left( \int_{(1-c)B_{n_k}^{-\frac{1}{\eta}}}^{B_{n_k}^{-\frac{1}{\eta}}} |W_{F_{\alpha,\eta}} (s,t)|^2 \frac{ds}{s} \right)^{\frac{p}{2}} dt \right]^{\frac{1}{p}} \nonumber \\
& \geq & K B_{n_k}^{-\alpha(1-2\ep_{j_k})+\frac{1}{\eta p}} \left[ \int_{X_{n_k}-cB_{n_k}}^{X_{n_k}+cB_{n_k}^{-\frac{1}{\eta}}} \left( \int_{(1-c)B_{n_k}^{-\frac{1}{\eta}}}^{B_{n_k}^{-\frac{1}{\eta}}} \frac{ds}{s} \right)^{\frac{p}{2}} dt \right]^{\frac{1}{p}} \nonumber \\
&=& K B_{n_k}^{-\alpha(1-2\ep_{j_k})+\frac{1}{\eta p}} \left[ \int_{X_{n_k}-cB_{n_k}}^{X_{n_k}+cB_{n_k}^{-\frac{1}{\eta}}} (-\ln(1-c))^{\frac{p}{2}} dt\right]^{\frac{1}{p}} \nonumber \\
&=& Kc \ln\left( \frac{1}{1-c}\right)^{1/2} B_{n_k}^{-\alpha(1-2\ep_{j_k})} \label{eqLFlb}
\EEN

\smallskip

\textbf{Conclusion : } Let $\ep \in \RR_+^*$. Contradiction, we suppose that $F\in T^{p}_{\frac{\alpha}{\delta}+\frac{1-\delta\eta}{\delta\eta p}+\ep} (x_0)$. According to the characterization supplied by Theorem \ref{propcar}, there exists $K'>0$ such that
\BE|L^{(p)}_{F_{\alpha,\eta}} (B_{n_k}^{-\frac{1}{\eta}},X_{n_k})| \leq K' B_{n_k}^{-\frac{1}{\eta}(\frac{\alpha}{\delta}+\frac{1-\delta\eta}{\delta\eta p}+\ep)}. \label{eqLFub}
\EE
According to the \eqref{eqLFlb} and \eqref{eqLFub},
\BE\label{ineqLp}
K B_{n_k}^{-\alpha(1-2\ep_{j_k})} \leq |L^{(p)}_{F_{\alpha,\eta}} (B_{n_k}^{-\frac{1}{\eta}},X_{n_k})| \leq K' B_{n_k}^{-\frac{1}{\eta}(\frac{\alpha}{\delta}+\frac{1-\delta\eta}{\delta\eta p}+\ep)}.
\EE
But $\delta\in [1,1/\eta)$, so
$$\frac{\alpha}{\delta}+\frac{1-\delta\eta}{\delta\eta p}+\ep \geq \alpha\eta+\ep \geq \alpha\eta$$
because $1-\delta\eta\in [0,1-\eta]$. But the last inequality implies that
$$\frac{1}{\eta}\left(\frac{\alpha}{\delta}+\frac{1-\delta\eta}{\delta\eta p}+\ep\right) >\alpha,$$
which contradicts \eqref{ineqLp} since $(\ep_j)_{j\in\NN^*}$ tends to $0$. Therefore, $F_{\alpha,\eta} \notin T^{(p)}_{\frac{\alpha}{\delta}+\frac{1-\delta\eta}{\delta\eta p}+\ep} (x_0)$ for all $\ep>0$ which proves that
$$h_{F_{\alpha,\eta}}^{(p)} (x_0) \leq \frac{\alpha}{\delta} + \frac{1-\delta\eta p}{\delta\eta p}.$$
\end{proof}

\subsection{Computation of the multifractal \texorpdfstring{$p$}{p}-spectrum}\label{subsec5.4}

According to Theorem 4.1 of \cite{Saes20} applied to the set $G'_1$ with this time $\rho>(3-3\alpha)/(1-\alpha\eta)$ and $\alpha<0$, we have almost surely that $G'_1$ is a covering of $[0,1]$. 
And so to conclude with the proof of the Theorem \ref{pthm1}, we just have to apply the propositions to determine inclusions on $p$-isohölderian sets. According to the Propositions \ref{pprop1} and \ref{pprop2}, almost surely, for all $x\in [0,1]$ and $\delta\in \left[1,\frac{1}{\eta}\right]$:
\begin{itemize}
\item[$\bullet$] If $h^{(p)}_{F_{\alpha,\eta}} (x) > \frac{\frac{\alpha\eta p +1}{\eta p}}{\delta}-\frac{1}{p}$ then $x\notin G'_\delta$.
\item[$\bullet$] If $h^{(p)}_{F_{\alpha,\eta}} (x) < \frac{\frac{\alpha\eta p +1}{\eta p}}{\delta}-\frac{1}{p}$ then $x\in G_\delta$.
\end{itemize}
Thus, almost surely for all $H\in\RR$
\BEN
h^{(p)}_{F_{\alpha,\eta}} (x)=H &\Longrightarrow & x \in \bigcap_{\delta< \frac{\frac{\alpha\eta p+1}{\eta p}}{H+\frac{1}{p}}} G_\delta \backslash \bigcup_{\delta > \frac{\frac{\alpha\eta p+1}{\eta p}}{H+\frac{1}{p}}} G'_\delta. \nonumber
\EEN
We deduce that almost surely, for all $H\in\RR$, 
\BE
    E_{F_{\alpha,\eta}} (H)\subset\bigcap_{\delta< \frac{\frac{\alpha\eta p+1}{\eta p}}{H+\frac{1}{p}}} G_\delta \backslash \bigcup_{\delta > \frac{\frac{\alpha\eta p+1}{\eta p}}{H+\frac{1}{p}}} G'_\delta. \nonumber
\EE
Conversely, according to the Propositions \ref{pprop1} and \ref{pprop2}, almost surely for any $H\in\RR$, 
\BE
    G'_{\frac{\alpha}{H}} \backslash \bigcup_{\delta > \frac{\frac{\alpha\eta p+1}{\eta p}}{H+\frac{1}{p}}} G_\delta \subseteq E_{F_{\alpha,\eta}} (H). \nonumber
\EE
For any $\delta\in [1,1/\eta]$, we have $H\in [\alpha\eta, \alpha+ \frac{1-\eta}{\eta p}]$ and thus we deduce for any $\ep>0$,
\BE\label{eq8.6p}
    D^{(p)}_{F_{\alpha,\eta}} (H)=\dim_H (E_{F_{\alpha,\eta}} (H)) \leq \dim_H \left(G_{\frac{\frac{\alpha\eta p+1}{\eta p}}{H+\frac{1}{p}}-\ep}\right) \leq \frac{1}{\frac{\frac{\alpha\eta p+1}{\eta p}}{H+\frac{1}{p}}-\ep} \nonumber
\EE
and
\BE
    D^{(p)}_{F_{\alpha,\eta}} (H) \geq \dim_H \left( G'_{\frac{\frac{\alpha\eta p+1}{\eta p}}{H+\frac{1}{p}}} \backslash \bigcup_{n\in\NN^*} G_{\frac{\frac{\alpha\eta p+1}{\eta p}}{H+\frac{1}{p}}+\frac{1}{n}} \right) \geq \frac{\eta p}{\alpha\eta p+1}\left(H+\frac{1}{p}\right). \nonumber
\EE
This allows us to conclude the point $(iii)$ of the Theorem \ref{pthm1}
$$D^{(p)}_F (H) = \frac{H\eta p+\eta}{\alpha\eta p +1}, \ H\in [\alpha\eta, \alpha+ \frac{1-\eta}{\eta p}].$$

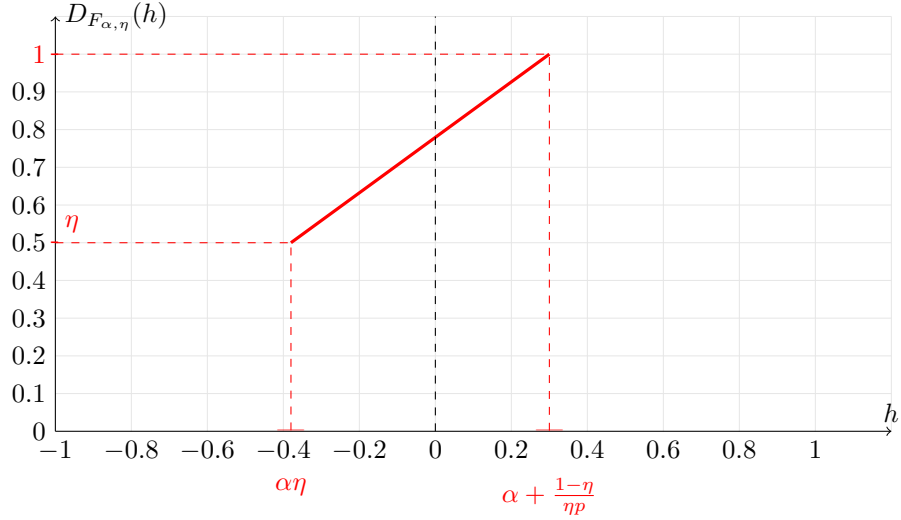
\begin{figure}
    \centering
    \begin{tikzpicture}[scale=0.5]

    \draw[-,very thin,gray!20] (2,0) -- (2,11);
    \draw[-,very thin,gray!20] (4,0) -- (4,11);
    \draw[-,very thin,gray!20] (6,0) -- (6,11);
    \draw[-,very thin,gray!20] (8,0) -- (8,11);
    \draw[-,very thin,gray!20] (10,0) -- (10,11);
    \draw[-,very thin,gray!20] (12,0) -- (12,11);
    \draw[-,very thin,gray!20] (14,0) -- (14,11);
    \draw[-,very thin,gray!20] (16,0) -- (16,11);
    \draw[-,very thin,gray!20] (18,0) -- (18,11);
    \draw[-,very thin,gray!20] (20,0) -- (20,11);
    \draw[-,very thin,gray!20] (22,0) -- (22,11);
    \draw[-,very thin,gray!20] (0,1) -- (22,1);
    \draw[-,very thin,gray!20] (0,2) -- (22,2);
    \draw[-,very thin,gray!20] (0,3) -- (22,3);
    \draw[-,very thin,gray!20] (0,4) -- (22,4);
    \draw[-,very thin,gray!20] (0,5) -- (22,5);
    \draw[-,very thin,gray!20] (0,6) -- (22,6);
    \draw[-,very thin,gray!20] (0,7) -- (22,7);
    \draw[-,very thin,gray!20] (0,8) -- (22,8);
    \draw[-,very thin,gray!20] (0,9) -- (22,9);
    \draw[-,very thin,gray!20] (0,10) -- (22,10);
    \draw[-,very thin,gray!20] (0,11) -- (22,11);
    \draw[->,black] (0,0) -- (22,0)  node[above] {$h$};
    \draw[->,black] (0,0) -- (0,11)  node[right] {$D_{F_{\alpha,\eta}} (h)$};
    
    \node[red] at (0,5) {-};
    \node[red] at (0,10) {-};
    \draw[red] (0,5.5) node[right] {$\eta$};
    \draw[red] (0,10) node[left] {$1$};
    \node[red] at (6.2,0) {|};
    \node[red] at (13,0) {|};
    \draw[red] (6.2,-1) node[below] {$\alpha\eta$};
    \draw[red] (13,-1) node[below] {$\alpha+\frac{1-\eta}{\eta p}$};
    
    \draw[-,red, very thick] (6.2,5) -- (13,10);

    \draw[dashed,red] (0,5) -- (6.2,5) -- (6.2,0);
    \draw[dashed,red] (0,10) -- (13,10) -- (13,0);
    \draw[dashed] (10,0) -- (10,11);

    \draw[] (0,0) node[below] {$-1$};
    \draw[] (2,0) node[below] {$-0.8$};
    \draw[] (4,0) node[below] {$-0.6$};
    \draw[] (6,0) node[below] {$-0.4$};
    \draw[] (8,0) node[below] {$-0.2$};
    \draw[] (10,0) node[below] {$0$};
    \draw[] (12,0) node[below] {$0.2$};
    \draw[] (14,0) node[below] {$0.4$};
    \draw[] (16,0) node[below] {$0.6$};
    \draw[] (18,0) node[below] {$0.8$};
    \draw[] (20,0) node[below] {$1$};

    \draw[] (0,0) node[left] {$0$};
    \draw[] (0,1) node[left] {$0.1$};
    \draw[] (0,2) node[left] {$0.2$};
    \draw[] (0,3) node[left] {$0.3$};
    \draw[] (0,4) node[left] {$0.4$};
    \draw[] (0,5) node[left] {$0.5$};
    \draw[] (0,6) node[left] {$0.6$};
    \draw[] (0,7) node[left] {$0.7$};
    \draw[] (0,8) node[left] {$0.8$};
    \draw[] (0,9) node[left] {$0.9$};

    \end{tikzpicture}
    \caption{Multifractal 1-spectrum of the sum of random pulses for$\alpha=-0.7$, $\eta=0.5$}
    \label{figpspec}
\end{figure}

According to the Definition \ref{defP} with the choice of the parameters $\eta\in (0,1)$, $\alpha\in (-\infty,0)$ such that $\eta-1<\alpha\eta$ and $p\in (1,-1/(\alpha\eta) + 1/\alpha\in (1/(\alpha\eta))$ in the Theorem \ref{pthm1}, we necessarily have the support of the multifractal spectrum which is a non empty segment containing $\alpha\eta <0$ and $\alpha+\frac{1-\eta}{\eta p}>0$. We represent on Figure \ref{figpspec}, an example of $p$-multifractal spectrum.

\section*{Conclusion}\label{sec6}

In conclusion, we presented a characterization of $p$-exponents by continuous wavelet transforms and then gave an example with the sums of random pulses. On the probabilistic side, a model extensively used and studied is supplied by Lévy processes  whose sample paths have locally bounded jumps and whose multifractal analysis has been performed in \cite{Jaffard99}. The multifractal spectrum of fractional integrals of Lévy processes were computed by   P. Balança in \cite{Balanca14}. A subsequent natural problem is then to determine  a  multifractal spectrum  of  fractional derivatives of Lévy processes. Fractional Lévy derivatives have been considered and studied in several papers \cite{Schilling01,Schilling97, Fageot17a,Fageot17b,Fageot20} and \cite{Unser14}. Since these processes have  non-locally bounded sample paths, in order to study their pointwise regularities, it is natural to use  continuous wavelets transforms.

\medskip

The study of continuous wavelet transform was also originally motivated by the study of the $p$-spectrum of Davenport series, which, given a sequence $(a_n)_{n\in\NN},$ is defined for every $x\in\mathbb{R}$ as 
\[D(x)=\sum_{n\in\NN} a_n \{nx\}\]
where 
\[ \{x\}= \left\{ \begin{tabular}{cl}
     $x-\lfloor x \rfloor-1/2$ &  if $x\notin \ZZ$ \\
     $0$ & otherwise.
\end{tabular}\right.\]

When $a_n=1/n^\beta$ with $\beta >2$, S. Jaffard established that if $s$ is a real in $(0,1)$ then the fractional derivative of order $s$ belongs to $L^p(\RR)$ as soon as $p<1/s$ and that the result is optimal in the sense that there are examples where $D$ is not in $L^p(\RR)$ when $p>1/s$. These fractional derivatives are an interesting example of non-locally bounded functions for which the multifractal analysis remains open.

\section*{Acknowledgments}\label{sec7}

The author thanks Stéphane Jaffard and Stéphane Seuret for all the enriching discussions around this article.

\bibliography{pleaders.bib}

%% BioMed_Central_Bib_Style_v1.01

\begin{thebibliography}{42}
% BibTex style file: bmc-mathphys.bst (version 2.1), 2014-07-24
\ifx \bisbn   \undefined \def \bisbn  #1{ISBN #1}\fi
\ifx \binits  \undefined \def \binits#1{#1}\fi
\ifx \bauthor  \undefined \def \bauthor#1{#1}\fi
\ifx \batitle  \undefined \def \batitle#1{#1}\fi
\ifx \bjtitle  \undefined \def \bjtitle#1{#1}\fi
\ifx \bvolume  \undefined \def \bvolume#1{\textbf{#1}}\fi
\ifx \byear  \undefined \def \byear#1{#1}\fi
\ifx \bissue  \undefined \def \bissue#1{#1}\fi
\ifx \bfpage  \undefined \def \bfpage#1{#1}\fi
\ifx \blpage  \undefined \def \blpage #1{#1}\fi
\ifx \burl  \undefined \def \burl#1{\textsf{#1}}\fi
\ifx \doiurl  \undefined \def \doiurl#1{\url{https://doi.org/#1}}\fi
\ifx \betal  \undefined \def \betal{\textit{et al.}}\fi
\ifx \binstitute  \undefined \def \binstitute#1{#1}\fi
\ifx \binstitutionaled  \undefined \def \binstitutionaled#1{#1}\fi
\ifx \bctitle  \undefined \def \bctitle#1{#1}\fi
\ifx \beditor  \undefined \def \beditor#1{#1}\fi
\ifx \bpublisher  \undefined \def \bpublisher#1{#1}\fi
\ifx \bbtitle  \undefined \def \bbtitle#1{#1}\fi
\ifx \bedition  \undefined \def \bedition#1{#1}\fi
\ifx \bseriesno  \undefined \def \bseriesno#1{#1}\fi
\ifx \blocation  \undefined \def \blocation#1{#1}\fi
\ifx \bsertitle  \undefined \def \bsertitle#1{#1}\fi
\ifx \bsnm \undefined \def \bsnm#1{#1}\fi
\ifx \bsuffix \undefined \def \bsuffix#1{#1}\fi
\ifx \bparticle \undefined \def \bparticle#1{#1}\fi
\ifx \barticle \undefined \def \barticle#1{#1}\fi
\bibcommenthead
\ifx \bconfdate \undefined \def \bconfdate #1{#1}\fi
\ifx \botherref \undefined \def \botherref #1{#1}\fi
\ifx \url \undefined \def \url#1{\textsf{#1}}\fi
\ifx \bchapter \undefined \def \bchapter#1{#1}\fi
\ifx \bbook \undefined \def \bbook#1{#1}\fi
\ifx \bcomment \undefined \def \bcomment#1{#1}\fi
\ifx \oauthor \undefined \def \oauthor#1{#1}\fi
\ifx \citeauthoryear \undefined \def \citeauthoryear#1{#1}\fi
\ifx \endbibitem  \undefined \def \endbibitem {}\fi
\ifx \bconflocation  \undefined \def \bconflocation#1{#1}\fi
\ifx \arxivurl  \undefined \def \arxivurl#1{\textsf{#1}}\fi
\csname PreBibitemsHook\endcsname

%%% 1
\bibitem[\protect\citeauthoryear{Frisch and Parisi}{1985}]{Frisch85}
\begin{bchapter}
\bauthor{\bsnm{Frisch}, \binits{U.}},
\bauthor{\bsnm{Parisi}, \binits{G.}}:
\bctitle{On the singularity structure of fully developed turbulence}.
In: \beditor{\bsnm{Ghil}, \binits{M.}},
\beditor{\bsnm{Benzi}, \binits{R.}},
\beditor{\bsnm{Parisi}, \binits{G.}} (eds.)
\bbtitle{Turbulence and Predictability in Geophysical Fluid Dynamics and
  Climate Dynamics},
pp. \bfpage{84}--\blpage{88}.
\bpublisher{North-Holland},
\blocation{New {Y}ork}
(\byear{1985})
\end{bchapter}
\endbibitem

%%% 2
\bibitem[\protect\citeauthoryear{Gagne}{1987}]{Gagne87}
\begin{botherref}
\oauthor{\bsnm{Gagne}, \binits{Y.}}:
Etude expérimentale de l’intermittence et des singularités dans le plan
  complexe en turbulence pleinement développée.
PhD thesis,
Université Joseph Fourier (France)
(1987)
\end{botherref}
\endbibitem

%%% 3
\bibitem[\protect\citeauthoryear{Calderón and Zygmund}{1961}]{Calderon61}
\begin{barticle}
\bauthor{\bsnm{Calderón}, \binits{A.}},
\bauthor{\bsnm{Zygmund}, \binits{A.}}:
\batitle{Local properties of solutions of elliptic partial differential
  equations}.
\bjtitle{Studia Mathematica}
\bvolume{20}(\bissue{2}),
\bfpage{181}--\blpage{225}
(\byear{1961})
\doiurl{10.1007/978-94-009-1045-4_17}
\end{barticle}
\endbibitem

%%% 4
\bibitem[\protect\citeauthoryear{Jaffard}{2006}]{Jaffard06a}
\begin{barticle}
\bauthor{\bsnm{Jaffard}, \binits{S.}}:
\batitle{Pointwise regularity associated with function spaces and multifractal
  analysis}.
\bjtitle{Banach Center Publications}
\bvolume{72},
\bfpage{93}--\blpage{100}
(\byear{2006})
\doiurl{10.4064/bc72-0-7}
\end{barticle}
\endbibitem

%%% 5
\bibitem[\protect\citeauthoryear{Jaffard}{2004}]{Jaffard04b}
\begin{barticle}
\bauthor{\bsnm{Jaffard}, \binits{S.}}:
\batitle{Wavelet techniques in multifractal analysis}.
\bjtitle{Proc. Symposia in Pure Mathematics}
\bvolume{72}(\bissue{2}),
\bfpage{91}--\blpage{152}
(\byear{2004})
\doiurl{10.1090/pspum/072.2/2112122}
\end{barticle}
\endbibitem

%%% 6
\bibitem[\protect\citeauthoryear{Holschneider and
  Tchamitchian}{1991}]{Holschneider91}
\begin{botherref}
\oauthor{\bsnm{Holschneider}, \binits{M.}},
\oauthor{\bsnm{Tchamitchian}, \binits{P.}}:
Pointwise analysis of {R}iemann's “nondifferentiable” function.
Inventiones mathematicae
\textbf{105}(1)
(1991)
\doiurl{10.1007/BF01232261}
\end{botherref}
\endbibitem

%%% 7
\bibitem[\protect\citeauthoryear{Jaffard}{1996}]{Jaffard96a}
\begin{barticle}
\bauthor{\bsnm{Jaffard}, \binits{S.}}:
\batitle{The spectrum of singularities of {R}iemann's function}.
\bjtitle{Revista Matemática Iberoamericana}
\bvolume{12}(\bissue{2}),
\bfpage{441}--\blpage{460}
(\byear{1996})
\doiurl{10.4171/RMI/203}
\end{barticle}
\endbibitem

%%% 8
\bibitem[\protect\citeauthoryear{Seuret and Ubis}{2017}]{Seuret17}
\begin{barticle}
\bauthor{\bsnm{Seuret}, \binits{S.}},
\bauthor{\bsnm{Ubis}, \binits{A.}}:
\batitle{Local $l^2$-regularity of {Riemann{\textquoteright}s} {Fourier}
  series}.
\bjtitle{Annales de l'Institut Fourier}
\bvolume{67}(\bissue{5}),
\bfpage{2237}--\blpage{2264}
(\byear{2017})
\doiurl{10.5802/aif.3135}
\end{barticle}
\endbibitem

%%% 9
\bibitem[\protect\citeauthoryear{Marmi et~al.}{1997}]{Marmi97a}
\begin{barticle}
\bauthor{\bsnm{Marmi}, \binits{S.}},
\bauthor{\bsnm{Moussa}, \binits{P.}},
\bauthor{\bsnm{Yoccoz}, \binits{J.-C.}}:
\batitle{The brjuno functions and their regularity properties}.
\bjtitle{Communications in Mathematical Physics}
\bvolume{186}(\bissue{2}),
\bfpage{256}--\blpage{293}
(\byear{1997})
\doiurl{10.1007/s002200050110}
\end{barticle}
\endbibitem

%%% 10
\bibitem[\protect\citeauthoryear{Jaffard and Martin}{2018}]{Jaffard18}
\begin{barticle}
\bauthor{\bsnm{Jaffard}, \binits{S.}},
\bauthor{\bsnm{Martin}, \binits{B.}}:
\batitle{Multifractal analysis of the {B}rjuno function}.
\bjtitle{Inventiones mathematicae}
\bvolume{212}(\bissue{1}),
\bfpage{109}--\blpage{132}
(\byear{2018})
\doiurl{10.1007/s00222-017-0763-z}
\end{barticle}
\endbibitem

%%% 11
\bibitem[\protect\citeauthoryear{Saës and Seuret}{Accepted}]{Saes20}
\begin{botherref}
\oauthor{\bsnm{Saës}, \binits{G.}},
\oauthor{\bsnm{Seuret}, \binits{S.}}:
Multifractal analysis of sums of random pulses.
Mathematical Proceedings of the Cambridge Philosophical Society
(Accepted)
\end{botherref}
\endbibitem

%%% 12
\bibitem[\protect\citeauthoryear{Arneodo et~al.}{2008}]{Arneodo08a}
\begin{botherref}
\oauthor{\bsnm{Arneodo}, \binits{A.}},
\oauthor{\bsnm{Audit}, \binits{B.}},
\oauthor{\bsnm{Kestener}, \binits{P.}},
\oauthor{\bsnm{Roux}, \binits{S.}}:
Wavelet-based multifractal analysis.
Scholarpedia
\textbf{3}(3)
(2008)
\doiurl{10.4249/scholarpedia.4103}
\end{botherref}
\endbibitem

%%% 13
\bibitem[\protect\citeauthoryear{Cioczek-Georges et~al.}{1995}]{Cioczek95a}
\begin{barticle}
\bauthor{\bsnm{Cioczek-Georges}, \binits{R.}},
\bauthor{\bsnm{Mandelbrot}, \binits{B.B.}},
\bauthor{\bsnm{Samorodnitsky}, \binits{G.}},
\bauthor{\bsnm{Taqqu}, \binits{M.S.}}:
\batitle{Stable fractal sums of pulses: the cylindrical case}.
\bjtitle{Bernoulli}
\bvolume{1}(\bissue{3}),
\bfpage{201}--\blpage{216}
(\byear{1995})
\doiurl{10.3150/bj/1193667815}
\end{barticle}
\endbibitem

%%% 14
\bibitem[\protect\citeauthoryear{Cioczek-Georges and
  Mandelbrot}{1995}]{Cioczek95b}
\begin{barticle}
\bauthor{\bsnm{Cioczek-Georges}, \binits{R.}},
\bauthor{\bsnm{Mandelbrot}, \binits{B.}}:
\batitle{A class of micropulses and antipersistent fractional brownian motion}.
\bjtitle{Stochastic Processes and their Applications}
\bvolume{60}(\bissue{1}),
\bfpage{1}--\blpage{18}
(\byear{1995})
\doiurl{10.1016/0304-4149(95)00046-1}
\end{barticle}
\endbibitem

%%% 15
\bibitem[\protect\citeauthoryear{Cioczek-Georges and
  Mandelbrot}{1996}]{Cioczek96a}
\begin{barticle}
\bauthor{\bsnm{Cioczek-Georges}, \binits{R.}},
\bauthor{\bsnm{Mandelbrot}, \binits{B.B.}}:
\batitle{Alternative micropulses and fractional brownian motion}.
\bjtitle{Stochastic Processes and their Applications}
\bvolume{64}(\bissue{2}),
\bfpage{143}--\blpage{152}
(\byear{1996})
\doiurl{10.1016/S0304-4149(96)00089-0}
\end{barticle}
\endbibitem

%%% 16
\bibitem[\protect\citeauthoryear{Lovejoy and Mandelbrot}{1985}]{Lovejoy85a}
\begin{barticle}
\bauthor{\bsnm{Lovejoy}, \binits{S.}},
\bauthor{\bsnm{Mandelbrot}, \binits{B.B.}}:
\batitle{Fractal properties of rain, and a fractal model}.
\bjtitle{Tellus A: Dynamic Meteorology and Oceanography}
\bvolume{37}(\bissue{3}),
\bfpage{209}--\blpage{232}
(\byear{1985})
\doiurl{10.3402/tellusa.v37i3.11668}
\end{barticle}
\endbibitem

%%% 17
\bibitem[\protect\citeauthoryear{Mandelbrot}{1995}]{Mandelbrot95a}
\begin{bchapter}
\bauthor{\bsnm{Mandelbrot}, \binits{B.B.}}:
\bctitle{Introduction to fractal sums of pulses}.
In: \beditor{\bsnm{Shlesinger}, \binits{M.F.}},
\beditor{\bsnm{Zaslavsky}, \binits{G.M.}},
\beditor{\bsnm{Frisch}, \binits{U.}} (eds.)
\bbtitle{L{\'e}vy Flights and Related Topics in Physics},
pp. \bfpage{110}--\blpage{123}.
\bpublisher{Springer},
\blocation{Berlin, Heidelberg}
(\byear{1995}).
\doiurl{10.1007/3-540-59222-9_29}
\end{bchapter}
\endbibitem

%%% 18
\bibitem[\protect\citeauthoryear{Demichel}{2006}]{Demichel06}
\begin{botherref}
\oauthor{\bsnm{Demichel}, \binits{Y.}}:
Analyse fractale d'une famille de fonctions al{\'e}atoires: les fonctions de
  bosses.
PhD thesis,
Universit{\'e} Blaise Pascal - Clermont-Ferrand II
(2006)
\end{botherref}
\endbibitem

%%% 19
\bibitem[\protect\citeauthoryear{Sa{\"e}s}{2021}]{Saes21}
\begin{botherref}
\oauthor{\bsnm{Sa{\"e}s}, \binits{G.}}:
Sommes fractales de pulses : {\'e}tude dimensionnelle et multifractale des
  trajectoires et simulations.
PhD thesis,
Universit{\'e} Paris-Est
(2021)
\end{botherref}
\endbibitem

%%% 20
\bibitem[\protect\citeauthoryear{Aubry and Jaffard}{2002}]{Aubry02}
\begin{barticle}
\bauthor{\bsnm{Aubry}, \binits{J.-M.}},
\bauthor{\bsnm{Jaffard}, \binits{S.}}:
\batitle{Random wavelet series}.
\bjtitle{Communications in Mathematical Physics}
\bvolume{227},
\bfpage{483}--\blpage{514}
(\byear{2002})
\doiurl{10.1007/s002200200630}
\end{barticle}
\endbibitem

%%% 21
\bibitem[\protect\citeauthoryear{Jaffard}{2000}]{Jaffard00}
\begin{barticle}
\bauthor{\bsnm{Jaffard}, \binits{S.}}:
\batitle{{On lacunary wavelet series}}.
\bjtitle{The Annals of Applied Probability}
\bvolume{10}(\bissue{1}),
\bfpage{313}--\blpage{329}
(\byear{2000})
\doiurl{10.1214/aoap/1019737675}
\end{barticle}
\endbibitem

%%% 22
\bibitem[\protect\citeauthoryear{Abry et~al.}{2015}]{Jaffard15a}
\begin{bchapter}
\bauthor{\bsnm{Abry}, \binits{P.}},
\bauthor{\bsnm{Jaffard}, \binits{S.}},
\bauthor{\bsnm{Leonarduzzi}, \binits{R.}},
\bauthor{\bsnm{Melot}, \binits{C.}},
\bauthor{\bsnm{Wendt}, \binits{H.}}:
\bctitle{Multifractal analysis based on $p$-exponents and lacunarity
  exponents}.
In: \beditor{\bsnm{Bandt}, \binits{C.}},
\beditor{\bsnm{Falconer}, \binits{K.}},
\beditor{\bsnm{Z{\"a}hle}, \binits{M.}} (eds.)
\bbtitle{Fractal Geometry and Stochastics V},
vol. \bseriesno{70},
pp. \bfpage{279}--\blpage{313}.
\bpublisher{Springer},
\blocation{Cham}
(\byear{2015}).
\doiurl{10.1007/978-3-319-18660-3_15}
\end{bchapter}
\endbibitem

%%% 23
\bibitem[\protect\citeauthoryear{Jaffard et~al.}{2016}]{Jaffard16}
\begin{barticle}
\bauthor{\bsnm{Jaffard}, \binits{S.}},
\bauthor{\bsnm{Melot}, \binits{C.}},
\bauthor{\bsnm{Leonarduzzi}, \binits{R.}},
\bauthor{\bsnm{Wendt}, \binits{H.}},
\bauthor{\bsnm{Abry}, \binits{P.}},
\bauthor{\bsnm{Roux}, \binits{S.G.}},
\bauthor{\bsnm{Torres}, \binits{M.E.}}:
\batitle{$p$-exponent and $p$-leaders, {P}art {I}: Negative pointwise
  regularity}.
\bjtitle{Physica A: Statistical Mechanics and its Applications}
\bvolume{448},
\bfpage{300}--\blpage{318}
(\byear{2016})
\doiurl{10.1016/j.physa.2015.12.061}
\end{barticle}
\endbibitem

%%% 24
\bibitem[\protect\citeauthoryear{Leonarduzzi et~al.}{2016}]{Leonarduzzi16}
\begin{barticle}
\bauthor{\bsnm{Leonarduzzi}, \binits{R.}},
\bauthor{\bsnm{Wendt}, \binits{H.}},
\bauthor{\bsnm{Abry}, \binits{P.}},
\bauthor{\bsnm{Jaffard}, \binits{S.}},
\bauthor{\bsnm{Melot}, \binits{C.}},
\bauthor{\bsnm{Roux}, \binits{S.G.}},
\bauthor{\bsnm{Torres}, \binits{M.E.}}:
\batitle{$p$-exponent and $p$-leaders, part ii: Multifractal analysis.
  relations to detrended fluctuation analysis}.
\bjtitle{Physica A: Statistical Mechanics and its Applications}
\bvolume{448},
\bfpage{319}--\blpage{339}
(\byear{2016})
\doiurl{10.1016/j.physa.2015.12.035}
\end{barticle}
\endbibitem

%%% 25
\bibitem[\protect\citeauthoryear{Perrier and Basdevant}{1996}]{Perrier96}
\begin{barticle}
\bauthor{\bsnm{Perrier}, \binits{V.}},
\bauthor{\bsnm{Basdevant}, \binits{C.}}:
\batitle{Besov norms in terms of the continous wavelet transform. application
  to structure functions}.
\bjtitle{Mathematical Models and Methods in Applied Sciences}
\bvolume{6}(\bissue{5}),
\bfpage{649}--\blpage{664}
(\byear{1996})
\doiurl{10.1142/S0218202596000262}
\end{barticle}
\endbibitem

%%% 26
\bibitem[\protect\citeauthoryear{Stein and Murphy}{1993}]{Stein93}
\begin{bbook}
\bauthor{\bsnm{Stein}, \binits{E.M.}},
\bauthor{\bsnm{Murphy}, \binits{T.S.}}:
\bbtitle{Harmonic Analysis : Real-Variable Methods, Orthogonality, and
  Oscillatory Integrals}.
\bpublisher{Princeton University Press},
\blocation{Princeton}
(\byear{1993})
\end{bbook}
\endbibitem

%%% 27
\bibitem[\protect\citeauthoryear{Jaffard and Mélot}{2005}]{Jaffard05}
\begin{barticle}
\bauthor{\bsnm{Jaffard}, \binits{S.}},
\bauthor{\bsnm{Mélot}, \binits{C.}}:
\batitle{Wavelet analysis of fractal boundaries. part 1: Local exponents}.
\bjtitle{Communications in Mathematical Physics}
\bvolume{258}(\bissue{3}),
\bfpage{513}--\blpage{565}
(\byear{2005})
\doiurl{10.1007/s00220-005-1354-1}
\end{barticle}
\endbibitem

%%% 28
\bibitem[\protect\citeauthoryear{Antoine et~al.}{2008}]{Murenzi08}
\begin{bbook}
\bauthor{\bsnm{Antoine}, \binits{J.-P.}},
\bauthor{\bsnm{Murenzi}, \binits{R.}},
\bauthor{\bsnm{Vandergheynst}, \binits{P.}},
\bauthor{\bsnm{Ali}, \binits{S.T.}}:
\bbtitle{Two-Dimensional Wavelets and Their Relatives}.
\bpublisher{Cambridge University Press},
\blocation{Cambridge}
(\byear{2008}).
\doiurl{10.1017/CBO9780511543395}
\end{bbook}
\endbibitem

%%% 29
\bibitem[\protect\citeauthoryear{Meyer}{1987}]{Meyer87}
\begin{bchapter}
\bauthor{\bsnm{Meyer}, \binits{Y.}}:
\bctitle{Principe d'incertitude, bases hilbertiennes et alg\`ebres
  d'op\'erateurs}.
In: \bbtitle{S\'eminaire Bourbaki : Volume 1985/86, Expos\'es 651-668}.
\bsertitle{Ast\'erisque},
pp. \bfpage{145}--\blpage{146}.
\bpublisher{Soci\'et\'e math\'ematique de France},
\blocation{Paris}
(\byear{1987})
\end{bchapter}
\endbibitem

%%% 30
\bibitem[\protect\citeauthoryear{Daubechies}{1992}]{Daubechies92}
\begin{bbook}
\bauthor{\bsnm{Daubechies}, \binits{I.}}:
\bbtitle{Ten Lectures on Wavelets}.
\bpublisher{Society for Industrial and Applied Mathematics},
\blocation{Philadelphie}
(\byear{1992}).
\doiurl{10.1137/1.9781611970104}
\end{bbook}
\endbibitem

%%% 31
\bibitem[\protect\citeauthoryear{Grossmann et~al.}{1990}]{Grossmann02}
\begin{bchapter}
\bauthor{\bsnm{Grossmann}, \binits{A.}},
\bauthor{\bsnm{Kronland-Martinet}, \binits{R.}},
\bauthor{\bsnm{Morlet}, \binits{J.}}:
\bctitle{Reading and understanding continuous wavelet transforms}.
In: \beditor{\bsnm{Combes}, \binits{J.-M.}},
\beditor{\bsnm{Grossmann}, \binits{A.}},
\beditor{\bsnm{Tchamitchian}, \binits{P.}} (eds.)
\bbtitle{Wavelets}.
\bpublisher{Springer},
\blocation{Berlin, Heidelberg}
(\byear{1990}).
\doiurl{10.1007/978-3-642-75988-8_1}
\end{bchapter}
\endbibitem

%%% 32
\bibitem[\protect\citeauthoryear{Jaffard et~al.}{1996}]{Jaffard96b}
\begin{bbook}
\bauthor{\bsnm{Jaffard}, \binits{S.}},
\bauthor{\bsnm{Meyer}, \binits{Y.}},
\bauthor{\bsnm{Ryan}, \binits{R.D.}}:
\bbtitle{Wavelet Methods for Pointwise Regularity and Local Oscillations of
  Functions}
vol. \bseriesno{123}.
\bpublisher{AMS},
\blocation{Providence}
(\byear{1996}).
\doiurl{10.1090/memo/0587}
\end{bbook}
\endbibitem

%%% 33
\bibitem[\protect\citeauthoryear{Meyer}{1990}]{Meyer90}
\begin{bbook}
\bauthor{\bsnm{Meyer}, \binits{Y.}}:
\bbtitle{Ondelettes et Opérateurs}.
\bpublisher{Hermann},
\blocation{Paris}
(\byear{1990})
\end{bbook}
\endbibitem

%%% 34
\bibitem[\protect\citeauthoryear{Jaffard}{1989}]{Jaffard89}
\begin{barticle}
\bauthor{\bsnm{Jaffard}, \binits{S.}}:
\batitle{Exposants de h{\"o}lder en des points donn{\'e}s et coefficients
  d’ondelettes}.
\bjtitle{CR Acad. Sci. Paris}
\bvolume{308}(\bissue{1}),
\bfpage{79}--\blpage{81}
(\byear{1989})
\end{barticle}
\endbibitem

%%% 35
\bibitem[\protect\citeauthoryear{Jaffard}{1999}]{Jaffard99}
\begin{barticle}
\bauthor{\bsnm{Jaffard}, \binits{S.}}:
\batitle{The multifractal nature of {L}évy processes}.
\bjtitle{Probability Theory and Related Fields}
\bvolume{114}(\bissue{2}),
\bfpage{207}--\blpage{227}
(\byear{1999})
\doiurl{10.1007/s004400050224}
\end{barticle}
\endbibitem

%%% 36
\bibitem[\protect\citeauthoryear{Balan{\c{c}}a}{2014}]{Balanca14}
\begin{barticle}
\bauthor{\bsnm{Balan{\c{c}}a}, \binits{P.}}:
\batitle{Fine regularity of {L}évy processes and linear (multi)fractional
  stable motion}.
\bjtitle{Electronic Journal of Probability}
\bvolume{19},
\bfpage{1}--\blpage{37}
(\byear{2014})
\doiurl{10.1214/EJP.v19-3393}
\end{barticle}
\endbibitem

%%% 37
\bibitem[\protect\citeauthoryear{Jacob and Schilling}{2001}]{Schilling01}
\begin{bbook}
\bauthor{\bsnm{Jacob}, \binits{N.}},
\bauthor{\bsnm{Schilling}, \binits{R.L.}}:
In: \beditor{\bsnm{Barndorff-Nielsen}, \binits{O.E.}},
\beditor{\bsnm{Resnick}, \binits{S.I.}},
\beditor{\bsnm{Mikosch}, \binits{T.}} (eds.)
\bbtitle{L{\'e}vy-Type Processes and Pseudodifferential Operators},
pp. \bfpage{139}--\blpage{168}.
\bpublisher{Birkh{\"a}user Boston},
\blocation{Boston, MA}
(\byear{2001}).
\doiurl{10.1007/978-1-4612-0197-7_7}
\end{bbook}
\endbibitem

%%% 38
\bibitem[\protect\citeauthoryear{Schilling}{1997}]{Schilling97}
\begin{barticle}
\bauthor{\bsnm{Schilling}, \binits{R.L.}}:
\batitle{On {F}eller processes with sample paths in {B}esov spaces}.
\bjtitle{Mathematische Annalen}
\bvolume{309}(\bissue{4}),
\bfpage{663}--\blpage{675}
(\byear{1997})
\doiurl{10.1007/s002080050132}
\end{barticle}
\endbibitem

%%% 39
\bibitem[\protect\citeauthoryear{Fageot et~al.}{2017a}]{Fageot17a}
\begin{barticle}
\bauthor{\bsnm{Fageot}, \binits{J.}},
\bauthor{\bsnm{Unser}, \binits{M.}},
\bauthor{\bsnm{Ward}, \binits{J.P.}}:
\batitle{On the {B}esov regularity of periodic {L}évy noises}.
\bjtitle{Applied and Computational Harmonic Analysis}
\bvolume{42}(\bissue{1}),
\bfpage{21}--\blpage{36}
(\byear{2017})
\doiurl{10.1016/j.acha.2015.07.001}
\end{barticle}
\endbibitem

%%% 40
\bibitem[\protect\citeauthoryear{Fageot et~al.}{2017b}]{Fageot17b}
\begin{barticle}
\bauthor{\bsnm{Fageot}, \binits{J.}},
\bauthor{\bsnm{Fallah}, \binits{A.}},
\bauthor{\bsnm{Unser}, \binits{M.}}:
\batitle{Multidimensional {L}évy white noise in weighted {B}esov spaces}.
\bjtitle{Stochastic Processes and their Applications}
\bvolume{127}(\bissue{5}),
\bfpage{1599}--\blpage{1621}
(\byear{2017})
\doiurl{10.1016/j.spa.2016.08.011}
\end{barticle}
\endbibitem

%%% 41
\bibitem[\protect\citeauthoryear{Aziznejad and Fageot}{2020}]{Fageot20}
\begin{barticle}
\bauthor{\bsnm{Aziznejad}, \binits{S.}},
\bauthor{\bsnm{Fageot}, \binits{J.}}:
\batitle{{Wavelet analysis of the {B}esov regularity of {L}évy white noise}}.
\bjtitle{Electronic Journal of Probability}
\bvolume{25},
\bfpage{1}--\blpage{38}
(\byear{2020})
\doiurl{10.1214/20-EJP554}
\end{barticle}
\endbibitem

%%% 42
\bibitem[\protect\citeauthoryear{Unser and Tafti}{2014}]{Unser14}
\begin{bbook}
\bauthor{\bsnm{Unser}, \binits{M.}},
\bauthor{\bsnm{Tafti}, \binits{P.D.}}:
\bbtitle{An Introduction to Sparse Stochastic Processes}.
\bpublisher{Cambridge University Press},
\blocation{Cambridge}
(\byear{2014}).
\doiurl{10.1017/CBO9781107415805}
\end{bbook}
\endbibitem

\end{thebibliography}

\end{document}